\pdfoutput=1

\documentclass[final]{siamltex}

\usepackage{graphicx}          		
\usepackage{graphics}         	 	
\usepackage{epsfig}            		
\usepackage{amsmath}           		
\usepackage{amssymb}           	
\usepackage[usenames,dvipsnames]{color}   
\usepackage[all]{xy}

\newtheorem{remark}[theorem]{remark}

\newcommand*\fvec[1]{\ensuremath{\mathbf{#1}}}  
\newcommand*\mc[0]{\mathcal}    
\newcommand*\mbb[0]{\mathbb}  
\newcommand{\until}[1]{\{1,\dots, #1\}}
\newcommand{\subscr}[2]{#1_{\textup{#2}}}

\newcommand*\terminalbus[0]{{{{\color{yellow}$\blacklozenge$}\hspace{-6pt}$\lozenge$\;}}}
\newcommand*\generator[0]{{\small{\color{red}$\blacksquare$}\hspace{-7.2pt}$\square$\;}}
\newcommand*\loadbus[0]{{\Large{\color{blue}$\bullet$}\hspace{-7.8pt}$\circ$\;}}

\newcommand\oprocendsymbol{\hbox{$\square$}}
\newcommand\oprocend{\relax\ifmmode\else\unskip\hfill\fi\oprocendsymbol}

\title{Kron Reduction of Graphs with Applications to Electrical Networks\thanks{This
    work was supported in part by NSF grants IIS-0904501 and CNS-0834446.}}

\title{Kron Reduction of Graphs\\ with Applications to Electrical Networks\thanks{This
    work was supported in part by NSF grants IIS-0904501 and CNS-0834446.}}

\author{Florian D\"orfler \and Francesco Bullo%
  \thanks{Florian D\"orfler and Francesco Bullo are with the Center for
    Control, Dynamical Systems and Computation, University of California at
    Santa Barbara, Santa Barbara, CA 93106, {\tt\small \{dorfler,
      bullo\}@engineering.ucsb.edu}} }

\begin{document}

\maketitle



\begin{abstract}
  Consider a weighted and undirected graph, possibly with self-loops, and
  its corresponding Laplacian matrix, possibly augmented with additional
  diagonal elements corresponding to the self-loops. The Kron reduction of
  this graph is again a graph whose Laplacian matrix is obtained by the Schur
  complement of the original Laplacian matrix with respect to a subset of
  nodes. 
  The Kron reduction process is ubiquitous in classic circuit theory
  and in related disciplines such as electrical impedance tomography, smart
  grid monitoring, transient stability assessment in power networks, or
  analysis and simulation of induction motors and power electronics. More
  general applications of Kron reduction occur in sparse matrix algorithms,
  multi-grid solvers, finite--element analysis, and Markov chains. The
  Schur complement of a Laplacian matrix and related concepts have 
  also been studied under different names and as purely theoretic problems in the literature on linear
  algebra.   
  In this paper we propose a general graph-theoretic framework for Kron reduction 
  that leads to novel and deep insights both on the mathematical and the physical side.  
  We show the applicability of our framework to various practical problem setups 
  arising in engineering applications and computation.
  Furthermore,
  we provide a comprehensive and detailed graph-theoretic analysis of the
  Kron reduction process encompassing topological, algebraic, spectral,
  resistive, and sensitivity analyses. Throughout our theoretic
  elaborations we especially emphasize the practical applicability of our
  results.
\end{abstract}




\pagestyle{myheadings}
\thispagestyle{plain}
\markboth{F. D\"orfler and F. Bullo}{Kron Reduction of Graphs with Applications to Electrical Networks}



\section{Introduction}\label{Section: Introduction}

Consider an undirected, connected, and weighted graph with $n$ nodes and adjacency matrix $A \in \mbb R^{n \times n}$. The corresponding loopy Laplacian matrix is the matrix $Q  \in \mathbb R^{n \times n}$ with off-diagonal elements $Q_{ij} = - A_{ij}$ and diagonal elements $Q_{ii} = A_{ii} + \sum_{j=1}^{n} A_{ij}$. Consider now a simple algebraic operation, namely the Schur complement of the loopy Laplacian matrix $Q$ with respect to a subset of nodes in the graph. As it turns out, the resulting lower dimensional matrix $\subscr{Q}{red}$ is again a well-defined loopy Laplacian matrix, and a graph can be naturally associated to it. 

This paper investigates this Schur complementation from the viewpoint of
algebraic graph theory. In particular we seek answers to the following
questions. How are the spectrum and the algebraic properties of $Q$ and
$\subscr{Q}{red}$ related? How about the corresponding graph topologies and
the effective resistances? What is the effect of a perturbation in the
original graph on the reduced graph, its loopy Laplacian $\subscr{Q}{red}$,
its spectrum, and its effective resistance? Finally, why is this graph
reduction process and its various properties of practical importance and in
which application areas? These are some of the questions that motivate this
paper.

\subsection*{Electrical Networks and the Kron Reduction}
To show the physical dimension and practical importance of the problem setup introduced above, we associate an electrical circuit to the graph induced by $A$. Consider a connected electrical network with $n$ nodes, branch conductances $A_{ij} \geq 0$, and shunt conductances $A_{ii} \geq 0$ connecting node $i$ to the ground. By Kirchhoff's and Ohm's laws the current-balance equations $I = Q V$ are obtained, where $I \in \mbb R^{n \times 1}$ are the currents injected at the nodes, $V \in \mbb R^{n \times 1}$ are the nodal voltages, and the conductance matrix $Q \in \mathbb R^{n \times n}$ is the loopy Laplacian matrix. In various applications of circuit theory and related disciplines it is desirable to obtain a lower dimensional\,\,electrically-equivalent network from the viewpoint of certain boundary nodes (or terminals) $\alpha \subsetneq \until n$, $|\alpha| \geq 2$. If $\beta = \until n \setminus \alpha$ denotes the set of interior nodes, then, after appropriately labeling the nodes, the current-balance equations can be partitioned as
\begin{equation}
	\left[ 	\begin{array}{c}
	I_{\alpha} \\ \hline I_{\beta}
	\end{array} \right]
	=
	\left[ 	\begin{array}{c|c}
	Q_{\alpha\alpha} & Q_{\alpha\beta}
	\\
	\hline
	Q_{\beta\alpha} & Q_{\beta\beta}
	\end{array} \right]
	\left[ 	\begin{array}{c}
	V_{\alpha} \\ \hline V_{\beta}
	\end{array} \right]
	\label{eq: network equations partitioned}
	\,.
\end{equation}
Gaussian elimination of the interior voltages $V_{\beta}$ in equations \eqref{eq: network equations partitioned} gives an electrically-equivalent reduced network with $|\alpha|$ nodes obeying the reduced current-balances
\begin{equation}
	I_{\alpha} + \subscr{Q}{ac}I_{\beta}
	=
	\subscr{Q}{red} V_{\alpha}
	\label{eq: reduced network equations}
	\,,
\end{equation}
where the reduced conductance matrix $\subscr{Q}{red} \in \mbb R^{|\alpha|\times|\alpha|}$ is again a loopy Laplacian given by the Schur complement of $Q$ with respect to the interior nodes $\beta$, that is, 
$
\subscr{Q}{red}
=
Q_{\alpha\alpha} - Q_{\alpha\beta} Q_{\beta\beta}^{-1} Q_{\beta\alpha}
$. The accompanying matrix $\subscr{Q}{ac} = -Q_{\alpha\beta}Q_{\beta\beta}^{-1} \in \mbb R^{|\alpha| \times (n-|\alpha|)}$ maps internal currents to boundary currents in the reduced network. In case that $I_{\beta}$ is the vector of zeros, the $(i,j)$-element of $\subscr{Q}{red}$ is the current at boundary node $i$ due to a unit potential at boundary node $j$ and a zero potential at all other boundary nodes. From here the reduced network can be further analyzed as an $|\alpha|$-port with current injections $I_{\alpha} + \subscr{Q}{ac}I_{\beta}$ and transfer conductance matrix $\subscr{Q}{red}$. 

\begin{figure}[htbp]
	\centering{
	\includegraphics[scale=0.3]{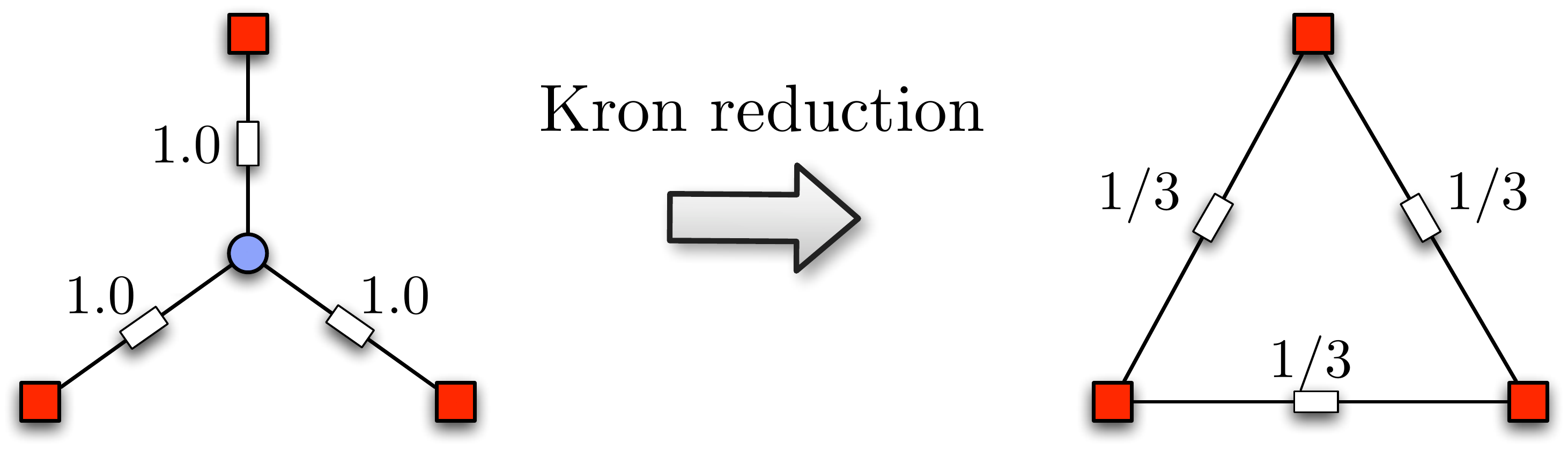}
	\caption{Kron reduction of a star-like electrical circuit with three boundary nodes \generator\!\!, one interior node \loadbus\!\!, and with unit conductances resulting in a reduced triangular reduced circuit.}
	\label{Fig: Y-Delta transform}
	}
\end{figure}

This reduction of an electrical network via a Schur complement of the
associated conductance matrix is known as {\it Kron reduction} due to the
seminal work of Gabriel Kron \cite{GK:39}, who identified fundamental
interconnections among physics, linear algebra, and graph theory
\cite{PMH:68,GK:53}.  The Kron reduction of a simple tree-like network
without current injections or shunt conductances is illustrated in Figure
\ref{Fig: Y-Delta transform}, an example familiar to every engineering
student as the $Y-\Delta$ transformation.

\subsection*{Literature Review}

The Kron reduction of networks is ubiquitous in circuit theory and related applications in order to obtain lower dimensional electrically-equivalent circuits. It appears for instance in the behavior, synthesis, and analysis of resistive circuits \cite{AJvdS:10,JCW-EIV:09,JCW:10}, particularly in the context of large-scale integration chips \cite{JR-WHAS:09,BNS:07,CSA-MHC-YII:05}. When applied to the impedance matrix of a circuit rather than the admittance matrix, Kron reduction is also referred to as the ``shortage operator'' \cite{WNA:71,WNA:75,PETJ-EJPP:09}. Kron reduction is a standard tool in the power systems community to obtain stationary and dynamically-equivalent reduced models for power flow studies \cite{JBW:09,HEB-RBS-DS-REN:07,AJW-BFW:96}, or in the reduction of differential-algebraic power network and RLC circuit models to lower dimensional purely dynamic models \cite{MAP:89,PWS-MAP:98,BA:02,FD-FB:10o-ieee,FD-FB:10o-ifac}. A recent application of Kron reduction is monitoring in smart power grids \cite{ID-MP:10} via synchronized phasor measurement units. Kron reduction is also crucial for reduced order modeling, analysis, and efficient simulation of induction motors \cite{RCD-MRG-SJS-DWB-RJN:02,OAM-GS-NA-SL-ZL:06} and power electronics \cite{KS-CE:06,AD:10}. Kron reduction is also known in the literature on electrical impedance tomography, where $\subscr{Q}{red}$ is referred to as the ``Dirichlet-to-Neumann map'' \cite{LB-VD-FGV:08,EC-M-JM:94,EC-DI-JM:98}. More generally, the Schur complement of a matrix and its associated graph is known in the context of (block) Gaussian elimination of sparse matrices \cite{MF:76,JRG:94,YS-MS:00,YS:03}, in sparse multi-grid solvers \cite{CW-WK-GW:97,AR:96}, and in finite-element analysis \cite{JM-MR-MT:00,RCD-MRG-SJS-DWB-RJN:02}. It serves as popular application example in linear algebra \cite{MF:86,YF:02,EAS-ARG:09,JJM-MN-HS-MJT:95}, a similar concept is employed in the stochastic complement \cite{CDM:89} or the cyclic reduction of Markov chains \cite{DAB-BM:09}, and a related concept is the Perron complement \cite{CDM:89-2,MN:00} of a matrix and its associated\,\,graph.

This brief literature review shows that Kron reduction is both a practically important and theoretically fascinating problem occurring in the reduction of networks and their associated matrices. Surprisingly, little is known about the graph-theoretic properties of the Kron reduction process. Yet the graph-theoretic analysis of the Kron reduction provides novel and deep insights both on the mathematical and the physical side of the considered problem. For instance, in the transient stability assessment of power networks, it is convenient to analyze the Kron-reduced network and afterwards relate its spectral and element-wise connectivity conditions to the non-reduced network \cite{FD-FB:10o-ieee,FD-FB:10o-ifac}. In electrical impedance tomography the Kron-reduced matrix $\subscr{Q}{red}$ can be directly obtained from boundary measurements or can be constructed from effective resistance measurements \cite{SN:02,AME:96}, and the problem is to infer properties of $Q$ from $\subscr{Q}{red}$ \cite{LB-VD-FGV:08,EC-M-JM:94,EC-DI-JM:98}. Other graph-theoretic problems occur in the applications of smart grid measurements \cite{ID-MP:10} and power flow studies \cite{JBW:09,HEB-RBS-DS-REN:07}, when disturbances such as load fluctuations or line outages in the non-reduced network have to be related to the Kron-reduced network or vice versa. Also the behavioral analysis of circuits relies strongly on Schur complement techniques and methods from graph theory \cite{AJvdS:10,JCW-EIV:09,JCW:10}. Yet another graph-theoretic problem occurring in numerous applications is minimizing the number of links during the reduction process to speed up costly computations associated with the network \cite{MF:76,JR-WHAS:09,YS-MS:00,YS:03,JRG:94,CW-WK-GW:97,DAB-BM:09,RCD-MRG-SJS-DWB-RJN:02,OAM-GS-NA-SL-ZL:06,AR:96,KS-CE:06,AD:10,JM-MR-MT:00}. 
In\,\,summary, it might be beneficial to analyze the Kron reduction process from a graph-theoretic viewpoint to address the problems above. Finally, the related literature on applied matrix analysis \cite{MF:86,YF:02,EAS-ARG:09,JJM-MN-HS-MJT:95,CDM:89,CDM:89-2,MN:00} shows that the graph-theoretic analysis of the Kron reduction process is an interesting mathematical problem in its own right.

\subsection*{Contributions}
The contributions of this paper are two-fold:

As a first contribution of this paper, we provide a graph-theoretic framework and interpretation of the Kron reduction process. This framework encompasses the applications of Kron reduction in classic circuit theory \cite{AJvdS:10,JCW-EIV:09,JCW:10,JR-WHAS:09,BNS:07,CSA-MHC-YII:05}, electrical impedance tomography \cite{LB-VD-FGV:08,SN:02,AME:96,EC-M-JM:94,EC-DI-JM:98}, power flow studies \cite{JBW:09,HEB-RBS-DS-REN:07,AJW-BFW:96}, monitoring in smart grids \cite{ID-MP:10}, and transient stability assessment \cite{MAP:89,PWS-MAP:98,BA:02,FD-FB:10o-ieee,FD-FB:10o-ifac}. We envision that our framework for Kron reduction captures various other reduction instances of physical problems, it can be applied to more abstract numerical problems as in \cite{MF:76,JRG:94,YS-MS:00,YS:03,CW-WK-GW:97,DAB-BM:09,AR:96}, and it can be extended to complex-valued, directed, and infinite-dimensional networks, possibly with an underlying dynamical system. 

The second contribution is a detailed analysis of the Kron reduction process. Essentially, Kron reduction of a connected graph, possibly with self-loops, is a Schur complement of corresponding loopy Laplacian matrix with respect to a subset of nodes. We relate the topological, the algebraic, and the spectral properties of the resulting Kron-reduced Laplacian matrix to those of the non-reduced Laplacian matrix. Furthermore, we relate the effective resistances in the original graph to the elements and effective resistances induced by the Kron-reduced Laplacian matrix. Thereby, we complement and extend various results in the literature on the effective resistance of a graph \cite{PETJ-EJPP:09,IG-WX:04,EB-AC-AME-JMG:09,PGD-JLS:84}. In our analysis, we carefully analyze the effects of self-loops, which model loads and dissipation in the applications of Kron reduction. Finally, we present a sensitivity analysis of the algebraic, spectral, and resistive properties of the Kron-reduced matrix with respect to perturbations in the non-reduced network topology. As a pure mathematical contribution, our analysis of Kron reduction complements the literature in linear algebra \cite{MF:86,YF:02,EAS-ARG:09,JJM-MN-HS-MJT:95} and analogous results on the Perron complement side \cite{CDM:89,CDM:89-2,MN:00}. In our analysis we do not aim at deriving only mathematical elegant results but also useful tools for the mentioned applications. We believe that our general analysis is a first step towards more detailed results in the specific application areas of Kron reduction.

\subsection*{Paper Organization} 
The remainder of this section introduces some notation recalls some preliminaries in matrix analysis, algebraic graph theory, and concerning the effective resistance. Section \ref{Section: Problem Setup and Applications} presents the general framework of Kron reduction and reviews various application areas and their specific problem setups. Section \ref{Section: Analysis of Kron Reduction} presents the graph-theoretic analysis of the Kron reduction process. Finally, Section \ref{Section: Conclusions} concludes the paper and suggests some future research directions.

\subsection*{Preliminaries and Notation}

{\em Sets:} 
Given a finite set $\mc Q$, let $|\mc Q|$ be its cardinality, and define for $n \in \mbb N$ the index set $\mc I_{n} = \until n$.

{\em Vectors and matrices:} 
Let $\fvec 1_{p \times q}$ and $\fvec 0_{p \times q}$ be the $p \times q$ dimensional matrices of unit and zero entries, and let $I_{n}$ be the $n$-dimensional identity matrix. For vectors, we adopt the shorthands $\fvec 1_{p} = \fvec 1_{p \times 1}$ and $\fvec 0_{p} = \fvec 0_{p \times 1}$ and define $e_{i}$ to be vector of zeros of appropriate dimension with entry $1$ at position $i$. For a real-valued 1d-array $\{x_{i}\}_{i=1}^{n}$, we let $\diag(\{x_{i}\}_{i=1}^{n}) \in \mbb R^{n \times n}$ be the associated diagonal matrix. 

Given a real-valued 2d-array $\{A_{ij}\}$ with $i,j \in \mc I_{n}$, let
$A \in \mbb R^{n \times n}$ denote the associated matrix and $A^{T}$ the transposed matrix. 
We use the following standard notation for submatrices \cite{FZ:05}: for two
non-empty index sets $\alpha,\beta \subseteq \mc I_{n}$ let
$A[\alpha,\beta]$ denote the submatrix of $A$ obtained by the rows indexed
by $\alpha$ and the columns indexed by $\beta$ and define the shorthands
$A[\alpha,\beta) = A[\alpha,\mc I_{n} \setminus \beta]$, $A(\alpha,\beta] =
A[\mc I_{n} \setminus \alpha,\beta]$, and $A(\alpha,\beta) = A[\mc I_{n}
\setminus \alpha,\mc I_{n} \setminus \beta]$.  We adopt the shorthand
$A[\{i\},\{j\}] = A[i,j] = A_{ij}$ for $i,j \in \mc I_{n}$, and for $x \in \mbb R^{n}$ the notation $x[\alpha,\{1\}] = x[\alpha]$ and $x(\alpha,\{1\}) = x(\alpha)$. For illustration, equation \eqref{eq: network equations partitioned} can be written unambiguously as
\begin{equation*}
	\left[ 	\begin{array}{c}
	I[\alpha] \\ \hline I(\alpha)
	\end{array} \right]
	=
	\left[ 	\begin{array}{c|c}
	Q[\alpha,\alpha] & Q[\alpha,\alpha)
	\\
	\hline
	Q(\alpha,\alpha] & Q(\alpha,\alpha)
	\end{array} \right]
	\left[ 	\begin{array}{c}
	V[\alpha] \\ \hline V(\alpha)
	\end{array} \right]
	\,.
\end{equation*}
If $A(\alpha,\alpha)$ is nonsingular, then the {\it Schur complement} of $A$ with respect to 
$A(\alpha,\alpha)$ (or equivalently the indices $\alpha$) is the
$|\alpha|\times|\alpha|$ dimensional matrix $A/A(\alpha,\alpha)$ defined by
\begin{equation*}
  A/A(\alpha,\alpha)
  \triangleq
  A[\alpha,\alpha] - A[\alpha,\alpha) A(\alpha,\alpha)^{-1} A(\alpha,\alpha] 
\,.
\end{equation*}
If $A$ is Hermitian, then we implicitly assume that its eigenvalues are arranged in increasing order: $\lambda_{1}(A) \!\leq\! \dots \!\leq\! \lambda_{n}(A)$. For a review of matrix analysis we refer to \cite{DSB:09}.

%

{\em Algebraic graph theory:} Consider the undirected, connected, and
weighted graph $G = (\mc I_{n},\mc E,A)$ with node set $\mc I_{n}$ and edge set $\mc E \subseteq \mc I_{n} \times \mc I_{n}$ induced by a symmetric,
nonnegative, and irreducible {\it adjacency matrix} $A \in \mbb R^{n \times
  n}$. A non-zero off-diagonal element $A_{ij} > 0$ corresponds to a weighted edge $\{i,j\} \in \mc E$, and a non-zero diagonal elements $A_{ii} > 0$ corresponds to a
weighted self-loop $\{i,i\} \in \mc E$.

The {\it Laplacian matrix} is the symmetric and irreducible matrix defined
by $L(A) = L \triangleq \diag(\{\sum_{j=1}^{n} A_{ij}\}_{i=1}^{n}) -
A$. Note that self-loops, even though apparent in the adjacency matrix $A$,
do not appear in the Laplacian matrix $L$. For these reasons\,\,and
motivated by the conductance matrix in circuit theory, we define the
symmetric and irreducible {\it loopy Laplacian matrix} $Q(A) = Q \triangleq
L + \diag(\{A_{ii}\}_{i=1}^{n}) \in \mbb R^{n \times n}$. Note that
adjacency matrix $A$ can be easily recovered from the loopy Laplacian $Q$
as $A = - Q + \diag(\{\sum_{j=1,j\neq i}^{n}Q_{ij}\}_{i=1}^{n})$, and thus
$Q$ uniquely induces the graph $G$. We refer to $Q$ as {\it strictly loopy}
Laplacian, respectively as {\it loop-less} Laplacian, if\,the\,\,graph
induced by $Q$ features at least one positively-weighted self-loop,
respectively no\,self-loops.

For a connected graph $\mathrm{ker}(L) = \mathrm{span}(\fvec 1_{n})$, and
all $n-1$ remaining non-zero eigenvalues of $L$ are strictly
positive. Specifically, the second-smallest eigenvalue $\lambda_{2}(L)$ is
a spectral connectivity measure called the {\it algebraic
  connectivity}. Recall that {\it irreducibility} of either $A$, $L$, or
$Q$ is equivalent to connectivity of $G$, which is again equivalent to
$\lambda_{2}(L) > 0$.
We refer the reader to \cite{CDG-GFR:01} for further details.

{\em Effective resistance:} The {\it effective resistance} 
$R_{ij}$ between two nodes $i,j \in \mc I_{n}$ of an undirected connected graph induced by a loopy Laplacian $Q$ is\,defined\,by
\begin{equation}
	R_{ij}
	\triangleq
	(e_{i}-e_{j})^{T} Q^{\dagger} (e_{i}-e_{j})
	=
	Q^{\dagger}_{ii} + Q^{\dagger}_{jj} - 2 Q^{\dagger}_{ij}
	\,,
	\label{eq: Definition of effective resistance}
\end{equation}
where $Q^{\dagger}$ is the Moore-Penrose pseudo inverse of $Q$. Since $Q^{\dagger}$ is symmetric (follows from the singular value decomposition), the matrix of effective resistances $R$ is again a symmetric matrix with zero diagonal elements $R_{ii}=0$. For two distinct nodes $i,j \in \mc I_{n}$ the reciprocal $1/R_{ij}$ is referred to as the {\it effective conductance} between\,\,$i$\,\,and\,\,$j$.

\begin{remark}[Physical interpretation]
\label{Remark: physical intuition about resistance}
\normalfont
If the graph is understood as a resistive electrical network with conductance matrix $Q$, the effective resistance $R_{ij}$ corresponds to the potential difference between the nodes $i$ and $j$ when a unit current is injected in $i$ and extracted in $j$. Indeed, in this case the current-balance equations  are $e_{i} - e_{j} = Q V$, where $V$ is the vector of potentials with respect to the ground (i.e., nodal voltages). The effective resistance $R_{ij}$ defined as the potential difference $R_{ij} = (e_{i} - e_{j})^{T} V$ can then be obtained via the impedance matrix $Q^{\dagger}$ as $R_{ij} = (e_{i}-e_{j})^{T} Q^{\dagger} (e_{i} - e_{j})$. 
\oprocend
\end{remark}

The effective resistance $R_{ij}$ can be thought of as a graph-theoretic distance between two nodes $i$ and $j$.  In most applications the effective resistance is defined for a pure topological network without shunt conductances, that is, a loop-less and uniformly weighted graph with $Q \equiv L$. We do not restrict ourselves to this case here. We refer the reader to \cite{PETJ-EJPP:09,FD-FB:10o-ifac,FF-AP-JMR-MS:07,IG-WX:04,AG-SB-AS:08,MF:98,EB-AC-AME-JMG:09} for various applications and properties of the effective resistance as well as interesting results relating $R$, $L$, $Q$, $L^{\dagger}$, and $Q^{-1}$.

\section{Problem Setup and Applications}\label{Section: Problem Setup and Applications}

\subsection{The Kron Reduction Process}\label{Subsection: Kron Reduction Process}

Consider an undirected, connected, and weighted graph $G = (\mc I_{n},\mc E,A)$ and its associated symmetric and irreducible matrices: the adjacency matrix $A \in \mbb R^{n \times n}$, Laplacian matrix $L(A)$, and loopy Laplacian matrix $Q(A)$. Furthermore, let $\alpha \subsetneq \mc I_{n}$ be a proper subset of nodes with $|\alpha| \geq 2$. We define the $(|\alpha| \times |\alpha|)$ dimensional {\it Kron-reduced matrix} $\subscr{Q}{red}$ by
\begin{equation}
	\subscr{Q}{red} \triangleq Q/Q(\alpha,\alpha) 
	= Q[\alpha,\alpha] - Q[\alpha,\alpha) Q(\alpha,\alpha)^{-1} Q(\alpha,\alpha] 
	\label{eq: Definition of Kron-reduced matrix}
	\,.
\end{equation}
In the following, we adopt the nomenclature from circuit theory and refer to the nodes  $\alpha$ and $\mc I_{n} \setminus \alpha$ as {\it boundary nodes} and {\it interior nodes}, respectively. The following lemma establishes the existence of the Kron-reduced matrix $\subscr{Q}{red}$ as well as the structural closure property of loopy Laplacian\,\,matrices.

\begin{lemma}[\bf Structural Properties of Kron Reduction]
\label{Lemma: Structural Properties of Kron Reduction}
Let $Q \in \mbb R^{n \times n}$ be a symmetric irreducible loopy Laplacian matrix and let $\alpha$ be a proper subset of $\mc I_{n}$ with $|\alpha| \geq 2$. 
Then the following statements hold:
\begin{enumerate}

	\item[1)] {\bf Existence:} The Kron-reduced matrix $\subscr{Q}{red} = Q/Q(\alpha,\alpha)$ is well defined.

	\item[2)] {\bf Closure:} If $Q$ is a symmetric loopy, strictly loopy, or loop-less Laplacian matrix, respectively, then $\subscr{Q}{red}$ is a symmetric loopy, strictly loopy, or loop-less Laplacian matrix, respectively.
	
	\item[3)] {\bf Accompanying matrix}: The accompanying matrix $\subscr{Q}{ac} \!\triangleq\! -Q[\alpha,\alpha)Q(\alpha,\alpha)^{-1}$\\$\in \mbb R^{|\alpha| \times (n-|\alpha|)}$ is nonnegative. If the subgraph among the interior nodes is connected and each boundary node is adjacent to at least one interior node, then $\subscr{Q}{ac}$ is positive. If additionally, $Q \equiv L$ is a loop-less Laplacian, then $\subscr{Q}{ac} = \subscr{L}{ac} \triangleq -L[\alpha,\alpha)L(\alpha,\alpha)^{-1}$ is column stochastic.

\end{enumerate}
\end{lemma}

\begin{figure}[htbp]
	\centering{
	\includegraphics[scale=0.215]{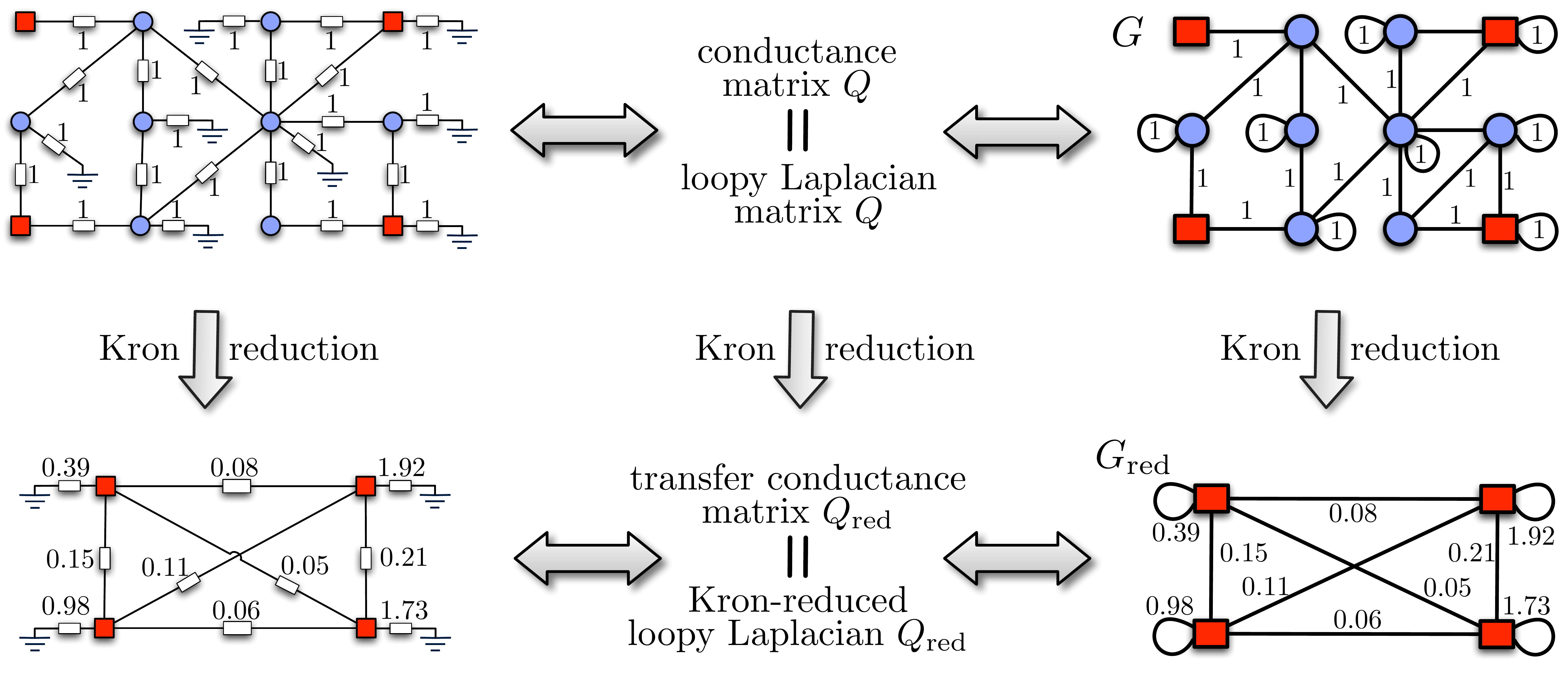}
	\caption{Illustration of an electrical network with 4 boundary nodes \generator\!\! and 8 interior nodes \loadbus\!\!, where all line and shunt conductances are of unit value. The associated conductance matrix (or loopy Laplacian) $Q$ and the corresponding graph $G$ are equivalent representations of this network. Kron reduction of the interior nodes \loadbus\! results in a reduced electrical network among the boundary nodes \generator\! with the Kron-reduced matrix $\subscr{Q}{red}$ as conductance matrix and the corresponding graph\,\,$\subscr{G}{red}$.}
	\label{Fig: Example for Kron reduction}
	}
\end{figure}

An interesting consequence of Lemma \ref{Lemma: Structural Properties of Kron Reduction} is that $\subscr{Q}{red}$, as a loopy Laplacian matrix, induces again an undirected and weighted graph. Hence, Kron reduction, originally defined as an algebraic operation in equation \eqref{eq: Definition of Kron-reduced matrix}, can be equivalently interpreted as a {\em graph-reduction process}, or as {\em physical reduction} of the associated\,\,circuit. This interplay between algebra, graph theory, and physics is illustrated in\,\,Figure\,\,\ref{Fig: Example for Kron reduction}.

In the following we denote the graph induced by $\subscr{Q}{red}$ as $\subscr{G}{red}$, and define $\subscr{A}{red} \triangleq - \subscr{Q}{red} + \diag(\{\sum_{j=1,j\neq i}^{n}\subscr{Q}{red}[i,j]\}_{i \in \alpha})$ and $\subscr{L}{red} \triangleq L(\subscr{A}{red})$ to be the corresponding reduced adjacency and loop-less Laplacian matrices. We remark that Lemma \ref{Lemma: Structural Properties of Kron Reduction} is partially noted in \cite{BA:02,MF:86,AJvdS:10,JCW-EIV:09,EAS-ARG:09,CW-WK-GW:97} and present its proof in the following.

{\em Proof of Lemma \ref{Lemma: Structural Properties of Kron Reduction}}.
By definition, $Q$ is (weakly) diagonally dominant since $Q_{ii} = \sum_{j=1,j \neq i}^{n} |Q_{ij}| + A_{ii}$ for all $i \in \mc I_{n}$. 
Due to the irreducibility of $Q$ the strict\,\,inequality 
$Q_{ii} 
> \sum_{j=1,j \neq i, j \not\in
  \alpha}^{n} |Q_{ij}| + A_{ii}$ holds at least for one $i \in \mc I_{n}
\setminus \alpha$. It follows that $Q(\alpha,\alpha)$ is also irreducible,
diagonally dominant, and\,\,has at least one row with strictly positive row
sum. Hence, $Q(\alpha,\alpha)$ is invertible
\cite[Corollary~6.2.27]{RAH-CRJ:85} and statement 1) follows.  

Statement 2)
is a consequence of the structural properties of $Q$ and the closure
properties of the Schur complement \cite[Chapter 4]{FZ:05}, which includes
the classes of symmetric, positive definite, and $M$-matrices. Since $Q$ is
a symmetric $M$-matrix, we directly conclude that $\subscr{Q}{red} =
Q/Q(\alpha,\alpha)$ is also a symmetric $M$-matrix. Hence,
$\subscr{Q}{red}$ is a symmetric loopy Laplacian matrix. This fact together
with the closure of positive definite matrices under the Schur complement
reveals that the class of symmetric strictly loopy Laplacian matrices is
closed under Kron reduction. To prove the closure of symmetric loop-less
Laplacians, assume without loss of generality that $\alpha=\mc
I_{|\alpha|}$, and consider the following equality for the row sums of the
loop-less Laplacian $Q$:
\begin{equation}
	\left[ 	\begin{array}{c|c}
	Q[\alpha,\alpha] & Q[\alpha,\alpha)
	\\
	\hline
	Q(\alpha,\alpha] & Q(\alpha,\alpha)
	\end{array} \right]
	\left[ 	\begin{array}{c}
	\fvec 1_{|\alpha|} \\ \hline \fvec 1_{n-|\alpha|}
	\end{array} \right]
	=
	\left[ 	\begin{array}{c}
	\fvec 0_{|\alpha|} \\ \hline \fvec 0_{|\alpha|}
	\end{array} \right]
	\label{eq: Equation determining row sums of Q}
	\,.
\end{equation}
Elimination of the second block of equations in \eqref{eq: Equation determining row sums of Q} results in $\fvec 0_{|\alpha|} = \subscr{Q}{red} \fvec 1_{|\alpha|}$, which shows that $\subscr{Q}{red}$ is a loop-less Laplacian and concludes the proof of statement 2). 

The second block of equations in \eqref{eq: Equation determining row sums of Q} can be rewritten as 
$\fvec 1_{n-|\alpha|} 
=  
\subscr{Q}{ac}^{T} \fvec 1_{|\alpha|} 
$, which shows that $\subscr{Q}{ac}$ is a column stochastic matrix in the loop-less case. In general, $\subscr{Q}{ac}= - Q[\alpha,\alpha) Q(\alpha,\alpha)^{-1}$ is nonnegative, since $-Q[\alpha,\alpha)$ and the inverse of the $M$-matrix $Q(\alpha,\alpha)$ are both nonnegative. If additionally each boundary node is connected to at least one interior node and the graph among the interior nodes is connected, then every row of $- Q[\alpha,\alpha)$ has at least one nonnegative entry and $Q(\alpha,\alpha)^{-1}$ is positive as well (since $Q(\alpha,\alpha)$ is an irreducible non-singular $M$-matrix \cite[Theorem 5.12]{MF:86}) , which guarantees positivity of $\subscr{Q}{ac}$. This completes the proof of statement 3).
\qquad\endproof

The Kron reduction process is ubiquitous in circuit theory and related scientific fields. Its general purpose is to construct low dimensional static or dynamic equivalents of higher dimensional models. In the following we describe different examples of Kron reduction in the fields of circuit theory, electrical impedance tomography, power flow studies, transient stability assessment, and smart grid monitoring.

\subsection{Kron Reduction in Large-Scale Integration Chips}

\begin{figure}[htbp]
	\centering{
	\includegraphics[scale=0.25]{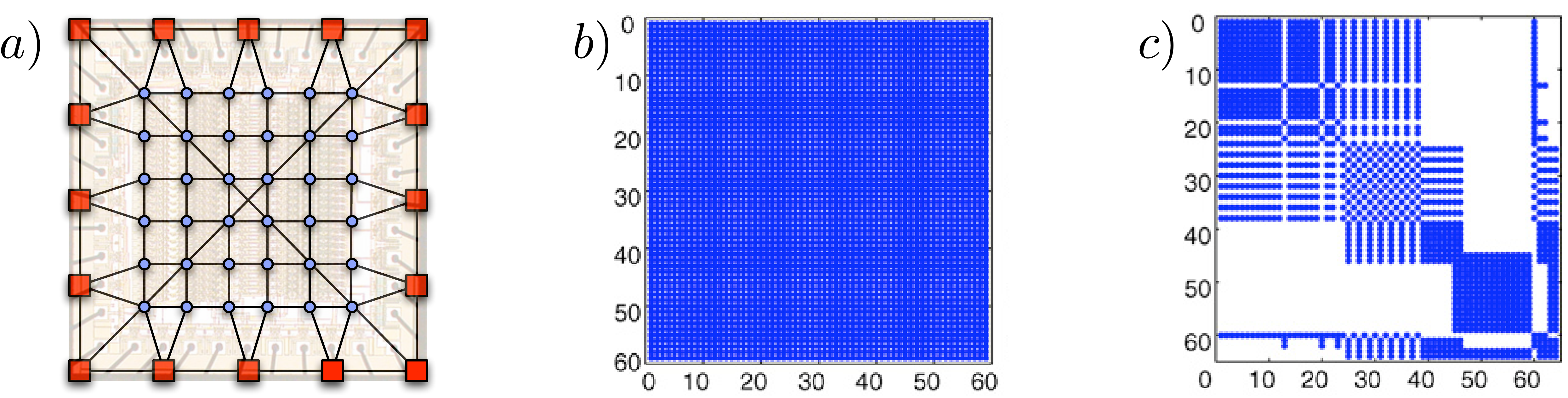}
	\caption{a) Illustration of an integration chip with a symmetric top-level network connecting the interior nodes with the terminals. The Kron reduction of all interior nodes of a network with 59 terminals results in a Kron-reduced matrix of dimension $59^{2}$ with $59^{2} = 3481$ non-zero entries. If all but five (specifically chosen) interior nodes are eliminated, then the Kron-reduced matrix is of dimension $64^{2}$ but has only $1506$ non-zero entries. The sparsity patterns of these two Kron-reduced matrices are illustrated in subfigures b) and c), which are taken from \cite{JR-WHAS:09}.}
	\label{Fig: integration chip}
	}
\end{figure}

In circuit theory, it is of interest to reduce the complexity of large-scale circuits by replacing them with an equivalent circuit with the same terminals (boundary nodes) but with fewer branches. This problem occurs in the context of large-scale integration chips \cite{JR-WHAS:09,BNS:07,CSA-MHC-YII:05}, where the smaller equivalent circuits are used to effectively compute the effective resistance between boundary nodes. The circuit reduction problem also stimulated a matrix-theoretic and behavioral analysis from the viewpoint of the boundary nodes \cite{AJvdS:10,JCW-EIV:09,JCW:10,BA:02} (where boundary nodes correspond to leaves in the modeling framework of \cite{JCW-EIV:09,JCW:10}).
For purely resistive circuits the proposed solution in \cite{AJvdS:10,JCW-EIV:09,JCW:10,JR-WHAS:09,BNS:07,CSA-MHC-YII:05} to obtain an equivalent circuit is the Kron reduction of the current-balance equations $I = Q V$ resulting in the reduced equations \eqref{eq: reduced network equations}. A particular reduction goal in \cite{JR-WHAS:09} is to reduce the fill-in of the Kron-reduced matrix $\subscr{Q}{red}$, see Figure \ref{Fig: integration chip}
for an illustration. The proper choice of the boundary nodes has a tremendous effect on the sparsity of the Kron-reduced matrix and saves numerical effort in subsequent computations, which is also a pervasive objective in the applications \cite{MF:76,JRG:94,YS-MS:00,YS:03,CW-WK-GW:97,JM-MR-MT:00,AR:96,RCD-MRG-SJS-DWB-RJN:02}. 

In \cite{JR-WHAS:09} a constructive algorithm is presented to obtain a sparse Kron-reduced matrix with minimal fill-in for computation of the effective resistance.
It is argued that reduction of a connected and sparse component of $Q$ results in a dense reduced component in $\subscr{Q}{red}$ and the effective resistance among boundary nodes is invariant under the Kron reduction process. We remark that these arguments are based on numerical observations and physical intuition rather than on mathematical proofs. This paper, puts the statements of \cite{JR-WHAS:09} on solid mathematical ground. We prove invariance of the effective resistance under Kron reduction and rigorously show how and under which conditions the topology changes from sparse to dense or even complete. The latter conditions may also be of interest for efficient matrix computations \cite{MF:76,JRG:94,YS-MS:00,YS:03,CW-WK-GW:97,JM-MR-MT:00,RCD-MRG-SJS-DWB-RJN:02,AR:96}. Moreover, our setup encompasses shunt loads  and currents drawn from the interior network, thereby generalizing results in \cite{AJvdS:10,JCW-EIV:09,JCW:10,BA:02}.

\subsection{Electrical Impedance Tomography}

\begin{figure}[htbp]
	\centering{
	\includegraphics[scale=0.535]{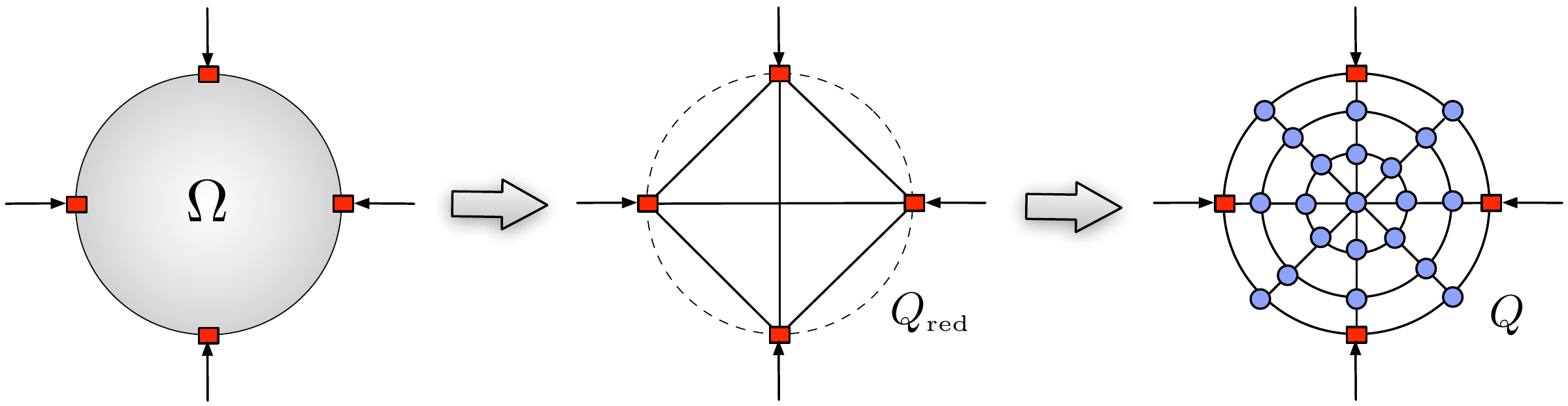}
	\caption{In electric impedance tomography the conductivity of the spatial domain $\Omega$ is estimated by measuring the Kron-reduced matrix $\subscr{Q}{red}$ at boundary ports \generator\!\!. From these measurements the conductance matrix $Q$ is re-constructed and serves as a spatial discretization of $\Omega$.}
	\label{Fig: Impedance tomography}
	}
\end{figure}

In electrical impedance tomography the goal is to determine the conductivity inside a compact and connected spatial domain $\Omega \subset \mathbb R^{2}$ from simultaneous measurements of currents and voltages at the boundary of $\Omega$, i.e., from measurement of the {\it Dirichlet-to-Neumann map}. Electrical impedance tomography finds applications in geophysics and medical imaging. A natural approach is a discretization of the spatial domain to a resistor network with conductance matrix $Q$ \cite{LB-VD-FGV:08}. When a unit potential is imposed at boundary node $j$ and a zero potential at all other boundary nodes, the current measured at boundary node $i$ gives $\subscr{Q}{red}[i,j]$, the $(i,j)$ element of the reduced conductance matrix (or discretized Dirichlet-to-Neumann map) $\subscr{Q}{red}$. Other methods iteratively construct the reduced impedance matrix (or discretized Neumann-to-Dirichlet map) $\subscr{Q}{red}^{\dagger}$ from measurements of the effective resistance matrix $R$ \cite{EC-M-JM:94,SN:02,AME:96}. The goal is then to recover the original network $Q$ from the reduced network $\subscr{Q}{red}$, i.e., inverting the Kron-reduction as in Figure \ref{Fig: Impedance tomography}. This is feasible only for highly symmetric networks as considered in \cite{EC-M-JM:94,EC-DI-JM:98,LB-VD-FGV:08}, but generally it is not possible to infer structural properties from $\subscr{Q}{red}$\,\,to\,\,$Q$.

This paper provides non-iterative methods to relate the effective resistance matrix $R$ and the inverse Kron-reduced matrix $\subscr{Q}{red}^{\dagger}$ as well as simple explicit formulas to relate $R$ and $\subscr{Q}{red}$ directly for uniform topologies. Furthermore, our analysis allows to partially invert the Kron reduction by estimating the spectrum of $Q$ or its effective resistance from the spectrum or resistance of $\subscr{Q}{red}$. Finally, our framework allows\,\,also for dissipation of energy in the spatial domain via discrete loads in the resistor network.

\subsection{Sensitivity of Reduced Power Flow}

\begin{figure}[htbp]
	\centering{
	\includegraphics[scale=0.119]{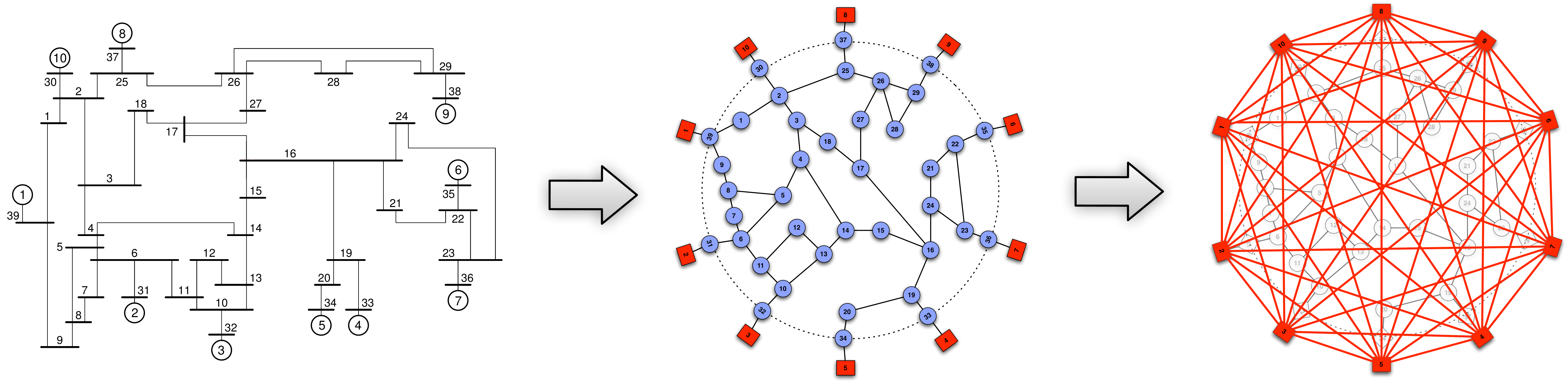}
	\caption{Single line diagram, of the New England Power Grid \cite{MAP:89}, an equivalent schematic representation with generators \generator\!\! and buses \loadbus\!\!, and the corresponding Kron-reduced network}
	\label{Fig: Reduction of New England Grid}
	}
\end{figure}

Large-scale power generation and transmission networks are also modeled as interconnected circuits \cite{PWS-MAP:98,AJW-BFW:96}. The nodes in a power grid are the 1-connected generators and the buses, see Figure \ref{Fig: Reduction of New England Grid}. 
Each edge connecting two nodes $i$ and $j$ is weighted by a non-zero line admittance $A_{ij}=A_{ji}\in\mbb C$ which is typically inductive, i.e., negative imaginary. Whereas a generator $i$ injects current $I_{i} \in \mbb C$ into the grid, the load at a bus $j$ draws a constant current $I_{j} \in \mbb C$ and features a constant shunt admittance $A_{jj} \in \mbb C$. Hence, the power network obeys the current-balance equations $I = Q V$, where the nodal {\it admittance matrix} $Q \in \mbb C^{n \times n}$ is the loopy Laplacian induced by the admittances $A_{ij}$. Alternatively, the current balance equations $I = Q V$ can be converted to the {\it power flow} or {\it power balance equations} $S = V \circ Q^{*}V^{*}$, where $\circ$ is the entry-wise (Hadamard) product, $^{*}$ denotes the conjugate transposed, and $S = V \circ I^{*}$ is the vector of power injections. Depending on the application and the set of available input data, it is more convenient to work with the linear current-balance equations or the nonlinear power-balance equations. 

A critical task in power network operation is monitoring, control, and optimization of the power flow to guarantee a stable and optimal operating point \cite{AJW-BFW:96}. Such an ideal operating point depends on the phases and magnitudes of the bus voltages and the powers and currents injected by the generators. The determining equations $S = V \circ Q^{*}V^{*}$ are too complicated to admit an analytic solution and often too onerous for a computational\,\,approach \cite{JBW:09,HEB-RBS-DS-REN:07}. If a set of boundary nodes $\alpha$ is identified for sensing or control purposes and all power injections are set to zero, then all interior nodes can be eliminated via Kron reduction leading to the reduced current-balance equations \eqref{eq: reduced network equations}.
The corresponding reduced power flow equations are obtained as $\subscr{S}{red} = V[\alpha] \circ \subscr{Q}{red}^{*}V[\alpha]^{*}$, where $\subscr{S}{red} = V[\alpha]\circ I[\alpha]^{*} + V[\alpha]\circ\subscr{Q}{ac}^{*}I(\alpha)^{*}$. 

For both analytic and computational purposes it is important to know the sensitivity of the reduced power flow w.r.t.\ changes in the original network topology and fluctuations in generation or load. For the lossless case when $Q$ is purely imaginary, this paper provides insightful and explicit results showing how perturbations in weights or topology of $Q$ affect the reduced {\it transfer admittance matrix} $\subscr{Q}{red}$. We also show the effect of shunt loads on the reduced network. In particular, a positive load $Q_{ii} > 0$ in the non-reduced network is shown to weaken the mutual transfer admittances $\subscr{Q}{red}[i,j]$ in the reduced network and to increase the reduced loads $\subscr{Q}{red}[i,i]$.

\subsection{Monitoring of DC Power Flow in Smart Grid}

\begin{figure}[htbp]
	\centering{
	\includegraphics[scale=0.72]{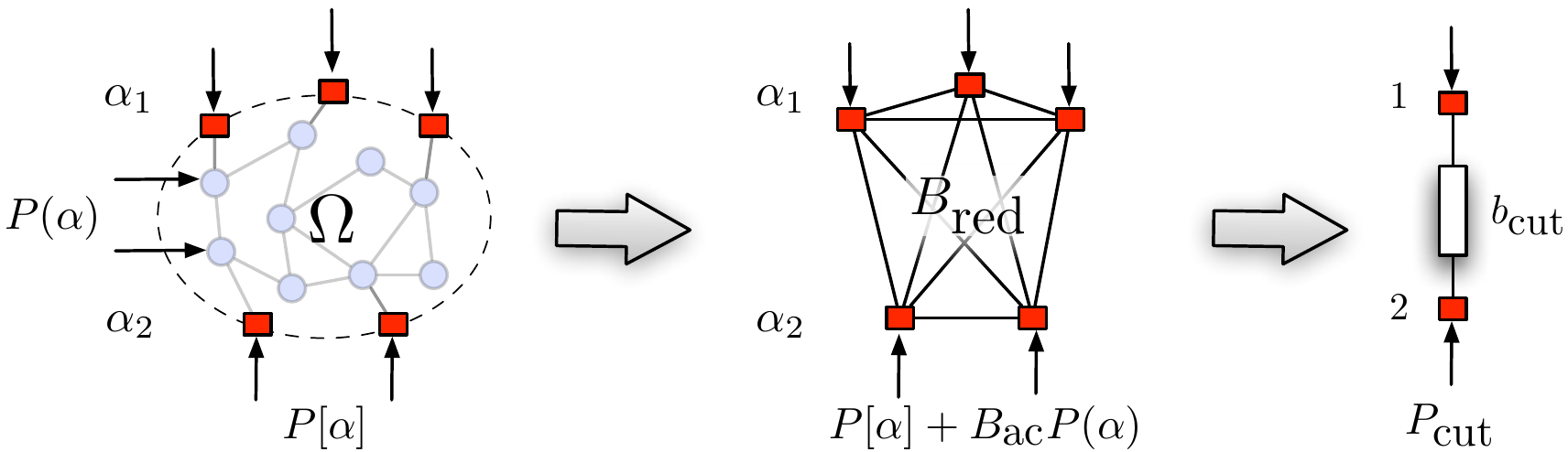}
	\caption{Reduction of the area $\Omega$ with boundary buses \generator\!\! to a single line equivalent $\{1,2\}$ describing the electrical characteristics between the boundary buses $\alpha_{1}$ and the boundary buses $\alpha_{2}$.}
	\label{Fig: cut set angle illustration}
	}
\end{figure}

The linearized {\it DC power flow equations} are $P = B \theta$, where $P = \Re(S) \in \mbb R^{n}$ is the vector of real power injections, $\theta \in \mbb R^{n}$ is the vector of voltage phase angles, and $B = - \Im(Q) \in \mbb R^{n \times n}$ is the {\it susceptance matrix}, i.e., the loopy Laplacian matrix induced by $-\Im(A_{ij})$. The DC power flow often serves as precursor for analytic and computational approaches to the nonlinear power flow \cite{AJW-BFW:96}. Consider now the problem of monitoring an area $\Omega$ of a {\it smart power grid} with $n$ nodes and equipped with synchronized phasor measurement units at the buses $\alpha = \{\alpha_{1},\alpha_{2}\}$  bordering the area \cite{ID-MP:10}. Kron reduction of the DC power flow $P = B \theta$ with respect to the interior nodes $\mc I_{n} \setminus \alpha$ yields the reduced DC power flow 
$P[\alpha] + \subscr{B}{ac} P(\alpha)= \subscr{B}{red} \theta[\alpha]$,
where $\subscr{P}{red}$ and $\subscr{P}{ac}$ are defined analogously to $\subscr{Q}{red}$ and $\subscr{Q}{ac}$. From here various scalar stress measures over the area $\Omega$ can be defined \cite{ID-MP:10}. Let $\sigma \in \mbb R^{|\alpha|}$ be the indicator vector for the boundary buses $\alpha_{1}$, that is, $\sigma_{i} = 1$ if $i \in \alpha_{1}$ and zero otherwise. The {\it cutset power flow} over the area $\Omega$ is $\subscr{P}{cut} \!=\! \sigma^{T} P[\alpha] + \sigma^{T} \subscr{B}{ac} P(\alpha)$, the {\it cutset susceptance} is $\subscr{b}{cut} \!=\! \sigma^{T} \subscr{B}{red} \sigma$, and the corresponding {\it cutset angle} is 
$\subscr{\theta}{cut} 
\!=\! 
\subscr{P}{cut}/\subscr{b}{cut} =
(\sigma^{T} P[\alpha] + \sigma^{T} \subscr{B}{ac} P(\alpha))/(\sigma^{T} \subscr{B}{red} \sigma)
$. Hence, the area $\Omega$ is effectively reduced to two nodes $\{1,2\}$ exchanging the power flow $\subscr{P}{cut}$ with angle $\subscr{\theta}{cut}$ over the susceptance $\subscr{b}{cut}$, see Figure \ref{Fig: cut set angle illustration}.
These scalar quantities indicate the stress within the area $\Omega$. For instance, a large cutset angle $\subscr{\theta}{cut}$ could be a blackout risk precursor. Of special interest are how load changes, line outages, or loss of entire nodes within the area $\Omega$ or on its boundary $\alpha$ affect the cutset angle\,\,$\subscr{\theta}{cut}$. 

This paper provides a comprehensive and detailed analysis of how changes in\,topol-ogy and weighting of the network affect the Kron-reduced matrix $\subscr{B}{red}$. These results include the self-loops in the graph (modeling loads) and can be easily translated to the cutset angle $\subscr{\theta}{cut}$ to show its sensitivity with respect to perturbations in the network.


\subsection{Transient Stability Assessment in Power Networks}
\label{Subsection: transient stability application}

Transient stability in a power network is the ability of the generators to remain in synchronism
when the network is subjected to large transient disturbances such as faults or loss of system components or severe fluctuations in generation or load. The dynamics of generator $i$ are typically given by the classic 
model \cite[Chapter 7]{PWS-MAP:98}
\begin{equation}
	M_{i} \ddot \theta_{i}
	=
	- D_{i} \dot{\theta}_{i} + P_{\textup{m},i} - P_{\textup{e},i} 
	,\quad
	i \in \mc I_{n}
	\label{eq: generator rotor dynamics}
	\,,
\end{equation}
where $(\theta_{i},\dot \theta_{i})$ are the generator rotor angle and speed, $P_{\textup{m},i} > 0$ is the mechanical power input, and $M_{i} > 0$ and $D_{i} > 0$ are the inertia and damping constant. The active output power injected by generator $i$ into the adjacent bus is $P_{\textup{e},i} = \Re(S_{i})$. The generator dynamics \eqref{eq: generator rotor dynamics} and the power flow equations $S = V \circ Q^{*}V^{*}$ define the classic structure-preserving power network model, an index-1 differential-algebraic system. We refer the reader to \cite{PWS-MAP:98,MAP:89,FD-FB:10o-ieee,FD-FB:10o-ifac} for further details and references. 

For the considered static load model, all bus nodes can be eliminated via Kron reduction and the network is reduced to the dynamic generator nodes. In the reduced network, the generators are coupled via the transfer admittance matrix $\subscr{Q}{red}$, and the electrical output power is $P_{\textup{e},i} = \Re(\subscr{S}{red}[i])$. For simplicity, we consider the lossless case $\Re(Q) = \fvec 0_{n \times n}$, which implies $\Re(\subscr{Q}{red}) = \fvec 0_{|\alpha| \times |\alpha|}$. Define the {\it effective power input} of each generator as $\omega_{i} = P_{\textup{m},i} + \Re(V_{i} \sum_{j=1}^{n-|\alpha|}\subscr{Q}{ac}^{*}[i,j] I_{|\alpha|+j}^{*}))$ and the {\it coupling weights}\,\,as $P_{ij} = |V_{i}| |V_{j}| \Im(\subscr{Q}{red}[i,j]) > 0$ (maximum power transfer between\,\,generators\,\,$i$\,\,and\,\,$j$), then equations \eqref{eq: generator rotor dynamics} together with the expression for $P_{\textup{e},i}$ give the {\it swing equations}
\begin{equation}
	M_{i} \ddot \theta_{i}
	=
	- D_{i} \dot{\theta}_{i} + \omega_{i}
	-
	\sum\nolimits_{j=1}^{|\alpha|} P_{ij} \sin(\theta_{i} - \theta_{j})
	\label{eq: Classic model - lossless}
	,\quad i \in \mc I_{n}
	\,.
\end{equation}
In \cite{FD-FB:09z}, we derived sufficient conditions under which the reduced model \eqref{eq: Classic model - lossless} {\it synchronizes}, i.e., all frequency differences $\dot \theta_{i}(t) - \dot \theta_{j}(t)$ converge to zero. For notational simplicity, we assume uniform damping here, that is, $D_{i} = D$ for all $i \in \alpha$. Then two sufficient conditions for synchronization of the model~\eqref{eq: Classic model - lossless} are \cite{FD-FB:09z}
\begin{equation}
 	|\alpha|  \min_{i \neq j} \{ P_{ij} \}
 	>\!
	\max_{i,j \in \mc I_{|\alpha|}} \left\{ \omega_{i}  - \omega_{j} \right\} \Bigg.
	 \,,
	\mbox{\;or \;\;}
  	\lambda_{2}(L(P_{ij}))
	>
	\Bigl( \sum\nolimits_{i,j=1,\,i<j}^{|\alpha|} (\omega_{i}-\omega_{j})^{2} \Bigr)^{\!1/2}
 	\label{eq: power network sync condition - reduced}
	\!\!.
\end{equation}
The left-hand sides of conditions \eqref{eq: power network sync condition - reduced} reflect the connectivity of the weighted graph induced by the\,\,power transfer $P_{ij}$: the term $|\alpha| \min_{i \neq j}\{ P_{ij} \} $ is a lower bound for $\min_{i} \sum_{j=1}^{|\alpha|} P_{ij}$, the worst coupling of one generator to the network, and $\lambda_{2}(L(P_{ij}))$ is the algebraic connectivity of the coupling. The right-hand side of conditions \eqref{eq: power network sync condition - reduced} describes the non-uniformity in effective power inputs $\omega_{i}$ in either $\infty$-norm or two-norm. In summary, conditions \eqref{eq: power network sync condition - reduced} can be interpreted as ``the reduced network connectivity has to dominate the network non-uniformity in effective power\,\,inputs.'' 

For uniformly lower-bounded voltages at all generators $|V_{i}| \geq V > 0$ (due to local control) the analysis of this paper will reveal that the spectral condition\,\,on $\lambda_{2}(L(P_{ij}))$ in the reduced network can be converted to the spectral synchronization condition
\begin{equation}
	\lambda_{2}(L)
	>
	\Bigl( \sum\nolimits_{i,j=1,\,i<j}^{|\alpha|} (\omega_{i}-\omega_{j})^{2} \Bigr)^{1/2} \frac{1}{V^{2}} + \max_{i \in \mc I_{n}}\{\subscr{A}{red}[i,i]\}
	\label{eq: power network spectral condition}
	\,,
\end{equation}
where $L$ is the Laplacian of the original lossless power network (weighted by $\Im(-A_{ij})$) and $\subscr{A}{red}[i,i]$ is the $i$th load in the reduced network. Similarly, if the effective resistance among all generators takes the uniform value $R$ and the effective resistance between the generators and the ground is uniform as well, then the results of this paper render the element-wise condition on $|\alpha| \min_{i \neq j} \{ P_{ij} \}$ in the reduced network to a resistive synchronization condition in the non-reduced\,\,network:
\begin{equation}
	\frac{1}{R}
 	>  
	\max_{i,j \in \mc I_{|\alpha|}} \left\{ \omega_{i}  - \omega_{j} \right\} \Bigg. \frac{1}{2V^{2}}
	+ \max_{i \in \mc I_{n}}\{\subscr{A}{red}[i,i]\}
	\label{eq: power network element-wise condition}
	\,.
\end{equation}
Conditions \eqref{eq: power network spectral condition} and \eqref{eq: power network element-wise condition} state that the network connectivity has to overcome the non-uniformity\,in\,effective power inputs and the dissipation by the loads, such that the network synchronizes.\,These conditions\,bridge the gap from classic transient stability analysis \cite{PWS-MAP:98,MAP:89} to the synchronization analysis in complex networks depending on the algebraic connectivity $\lambda_{2}(L)$ \cite{AA-ADG-JK-YM-CZ:08} or the effective conductance $1/R$\,\,\cite{GK-MBH-KEB-MJB-BK-DA:06}.

\section{Kron Reduction of Graphs}\label{Section: Analysis of Kron Reduction}

This section analyzes the algebraic, topological, spectral, and sensitivity properties, as well as the effective resistance of the Kron-reduced matrix $\subscr{Q}{red}$ and its associated graph. 
Throughout this section we  assume that $Q \in \mbb R^{n \times n}$ is a symmetric and irreducible loopy Laplacian matrix (corresponding to an undirected and connected graph with $n$ nodes), and we let $\alpha$ be a proper subset of $\mc I_{n}$ with $|\alpha| \geq 2$. For simplicity and without loss of generality, we assume that the $n$ nodes are labeled such that $\alpha = \mc I_{|\alpha|}$.
All results in this paper can also be stated for an arbitrary labeling of the nodes at the cost of more complicated notation.

\subsection{Preliminary Results on the Augmented Laplacian Matrix and Iterative Kron Reduction}

The following definitions of the {\it augmented Laplacian} matrix and the {\it iterative Kron reduction} will be central to the subsequent developments both for illustration and analysis purposes.

\begin{figure}[htbp]
	\centering{
	\includegraphics[scale=0.24]{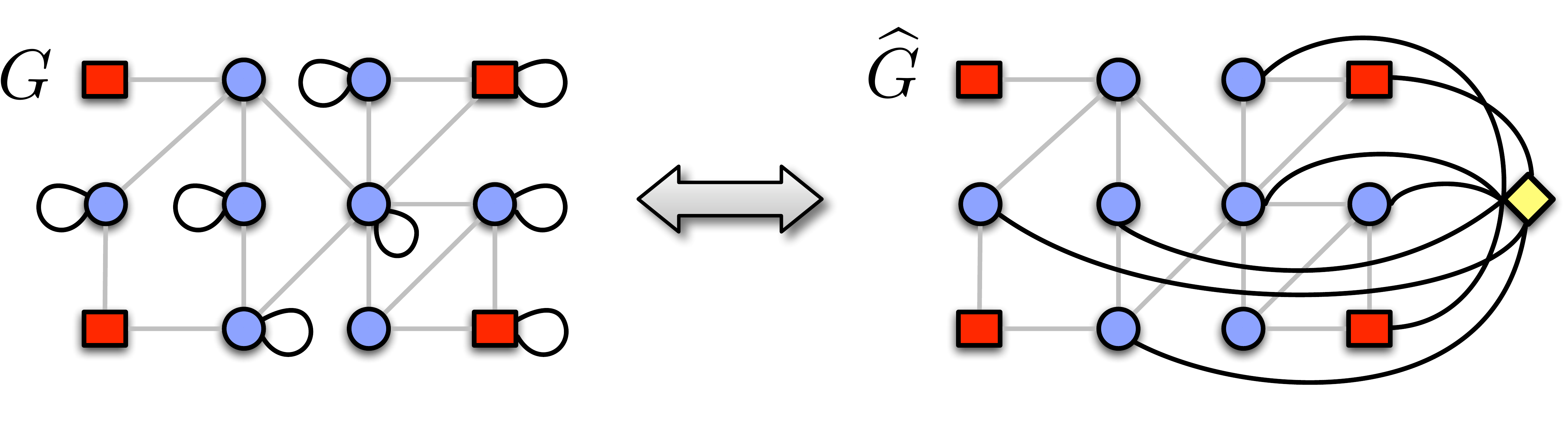}
	\caption{Illustration of the graph $G$ associated with the electrical network from Figure \ref{Fig: Example for Kron reduction} 
	and the corresponding augmented graph $\widehat G$ with additional grounded node \terminalbus\!\!.}
	\label{Fig: Example for loopy Kron reduction}
	}
\end{figure}

The role of the self-loops induced by a {\em strictly} loopy Laplacian $Q \in \mbb R^{n \times n}$ can be better understood by introducing the {\it grounded node} with index $n+1$. Then the strictly loopy Laplacian $Q$ can be regarded as the principal $n\times n$ block embedded in the $(n+1) \times (n+1)$ dimensional {\it augmented Laplacian matrix}
\begin{equation}
	\widehat Q
	\triangleq
	\left[\begin{array}{c|c}
	Q & - \diag(\{A_{ii}\}_{i=1}^{n}) \fvec 1_{n}
	\\\hline
	- \fvec 1_{n}^{T} \diag(\{A_{ii}\}_{i=1}^{n}) & \sum_{i=1}^{n} A_{ii}
	\end{array}\right]
	\label{eq: augmented Laplacian}
	\,,
\end{equation}
where $A \in \mbb R^{n \times n}$ is the adjacency matrix corresponding to $Q$. The augmented Laplacian $\widehat Q$ is the Laplacian of the augmented graph $\widehat G$ with node set $\widehat {\mc V} = \{\mc I_{n},n+1\}$ and edge set $\widehat {\mc E} = \{\mc E,\subscr{\mc E}{augment}\}$.
Here a node $i \in \mc I_{n}$ is connected to the grounded node $n+1$ via a weighted edge $\{i,n+1\} \in \subscr{\mc E}{augment}$ if and only if $A_{ii} > 0$. A graph $G$ and its associated augmented graph $\widehat G$ are  illustrated in Figure \ref{Fig: Example for loopy Kron reduction}.

This augmentation and embedding process is\,inspired by Fiedler's work \cite{MF:98}\,on\,non-singular $M$-matrices.
The following lemma presents some interesting properties of $\widehat Q$.

\begin{lemma}[\bf Properties of augmented Laplacian]
\label{Lemma: Properties of augmented Laplacian matrix}
Consider the symmetric and irreducible strictly loopy Laplacian $Q \in \mbb R^{n \times n}$ and the corresponding augmented Laplacian matrix $\widehat Q \in \mbb
R^{(n+1) \times (n+1)}$ defined in \eqref{eq: augmented Laplacian}. 
The following statements hold:
\begin{enumerate}

\item[1)] {\bf Algebraic properties:} $\widehat Q$ is an irreducible and
  symmetric loop-less Laplacian matrix.

\item [2)] {\bf Spectral properties:} The eigenvalues of $Q$ and $\widehat
  Q$ interlace each other:
  \begin{equation*}
    0 = \lambda_{1}(\widehat Q) < \lambda_{1}(Q) \leq \lambda_{2}(\widehat Q) \leq \lambda_{2}(Q) \leq \dots \leq \lambda_{n}(\widehat Q) \leq \lambda_{n}(Q) \leq \lambda_{n+1}(\widehat Q)
    \,.
  \end{equation*}
  
\item[3)] {\bf Kron reduction:} 
  Consider the strictly loopy Laplacian $\subscr{Q}{red} \!=\! Q/Q(\alpha,\alpha)$\,and the loop-less
  Laplacian $\subscr{\widehat Q}{red} \!\triangleq\! \widehat Q/\widehat
  Q(\{\alpha,n+1\}\!,\!\{\alpha,n+1\})$, both obtained\,by\,\,Kron reduction of the
  interior nodes $\mc I_{n} \setminus \alpha$. Then the following diagram commutes:
	\begin{center}
		\includegraphics[scale=0.32]{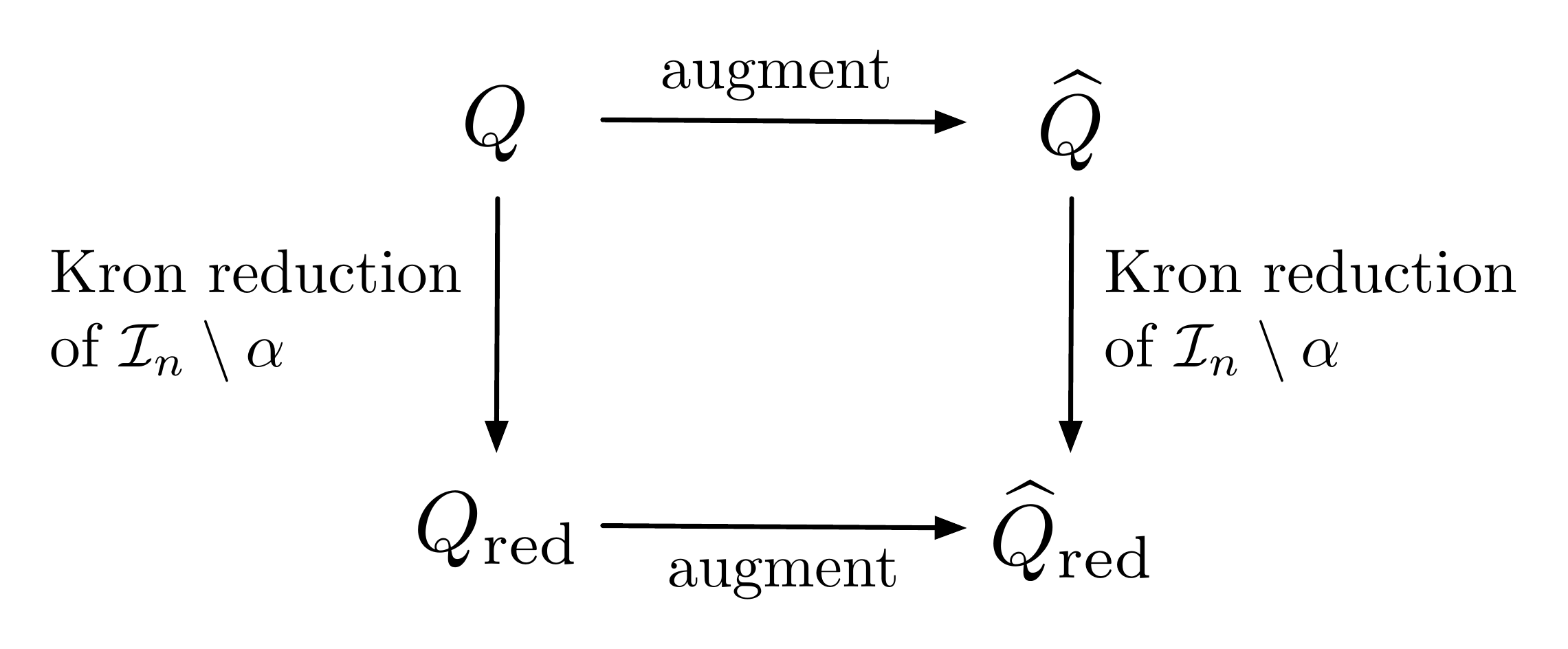}
	\end{center}	
	In equivalent words, $\subscr{\widehat Q}{red}$ is the augmented Laplacian 
        associated to $\subscr{Q}{red}$, i.e.,
	\begin{equation}
	\subscr{\widehat Q}{red}
	\equiv
	\left[\begin{array}{c|c}
	\subscr{Q}{red} & - \diag(\{\subscr{A}{red}[i,i]\}_{i \in \alpha}) \fvec 1
	\\\hline
	- \fvec 1^{T} \diag(\{\subscr{A}{red}[i,i]\}_{i \in \alpha}) & \sum_{i=1}^{n} \subscr{A}{red}[i,i]
	\end{array}\right]
	\label{eq: augmented Laplacian reduced}
	\,.
	\end{equation}
\end{enumerate}
\end{lemma}

Properties 2) and 3) of Lemma \ref{Lemma: Properties of augmented Laplacian matrix} intuitively illustrate the effect of self-loops on the spectrum of $Q$ and its Kron-reduced matrix. Specifically, the elegant relationship 3) implies that the reduced self-loops $\subscr{A}{red}[i,i]$ can be equivalently obtained by the last column or row of $\subscr{\widehat Q}{red}$. In short, Kron reduction can equivalently be applied to the strictly loopy network $G$ or to the augmented loop-less network $\widehat G$.

{\em Proof of Lemma \ref{Lemma: Properties of augmented Laplacian matrix}}.
Property 1) follows trivially from the construction of the augmented Laplacian $\widehat Q$. Property 2) is a direct application of the {\it interlacing theorem for bordered matrices} \cite[Theorem 4.3.8]{RAH-CRJ:85}, where $0 = \lambda_{1}(\widehat Q) < \lambda_{1}(Q)$ since $\widehat Q$ is an irreducible loop-less Laplacian and $Q$ is non-singular. In property 3), the (1,1) block of the matrix on the right-hand side of identity \eqref{eq: augmented Laplacian reduced} follows by writing out the Schur complement of a matrix partitioned in $3 \times 3$ blocks, as in the proof of the {\it Quotient Formula} \cite[Theorem 1.4]{FZ:05}. The other blocks follow immediately since Kron reduction of the loop-less Laplacian $\subscr{\widehat Q}{red}$ yields again a loop-less Laplacian by Lemma \ref{Lemma: Structural Properties of Kron Reduction}. 
\qquad\endproof

Gaussian elimination of interior voltages from the current-balance equations $I = Q V$ can either be performed via Kron reduction in a single step, as in equation \eqref{eq: reduced network equations}, or in multiple steps, each interior node $\ell \in \until{n-|\alpha|}$ at a time. The following concept of {\it iterative Kron reduction} addresses exactly this point.

\begin{definition}[\bf Iterative Kron reduction]
\label{Definition: iterative Kron reduction}
{\it Iterative Kron reduction} associates to a symmetric irreducible loopy
Laplacian matrix $Q \in \mbb R^{n \times n}$ and indices
$\until{|\alpha|}$, a sequence of matrices $Q^{\ell} \in \mbb R^{(n - \ell)
  \times (n - \ell)}$, $\ell\in\{1,\dots,n-|\alpha|\}$,\,\,defined\,\,by
\begin{equation}
	Q^{\ell}
	= 
	Q^{\ell-1}/Q^{\ell-1}_{k_{\ell} k_{\ell}}
	\label{eq: iterative Kron reduction in vector form}
	\,,
\end{equation}
where $Q^{0} = Q$ and $k_{\ell} = n+1-\ell$, that is, $Q^{\ell-1}_{k_{\ell} k_{\ell}}$ is the lowest diagonal entry of $Q^{\ell-1}$.
\end{definition}

\begin{figure}[htbp]
	\centering{
	\includegraphics[scale=0.182]{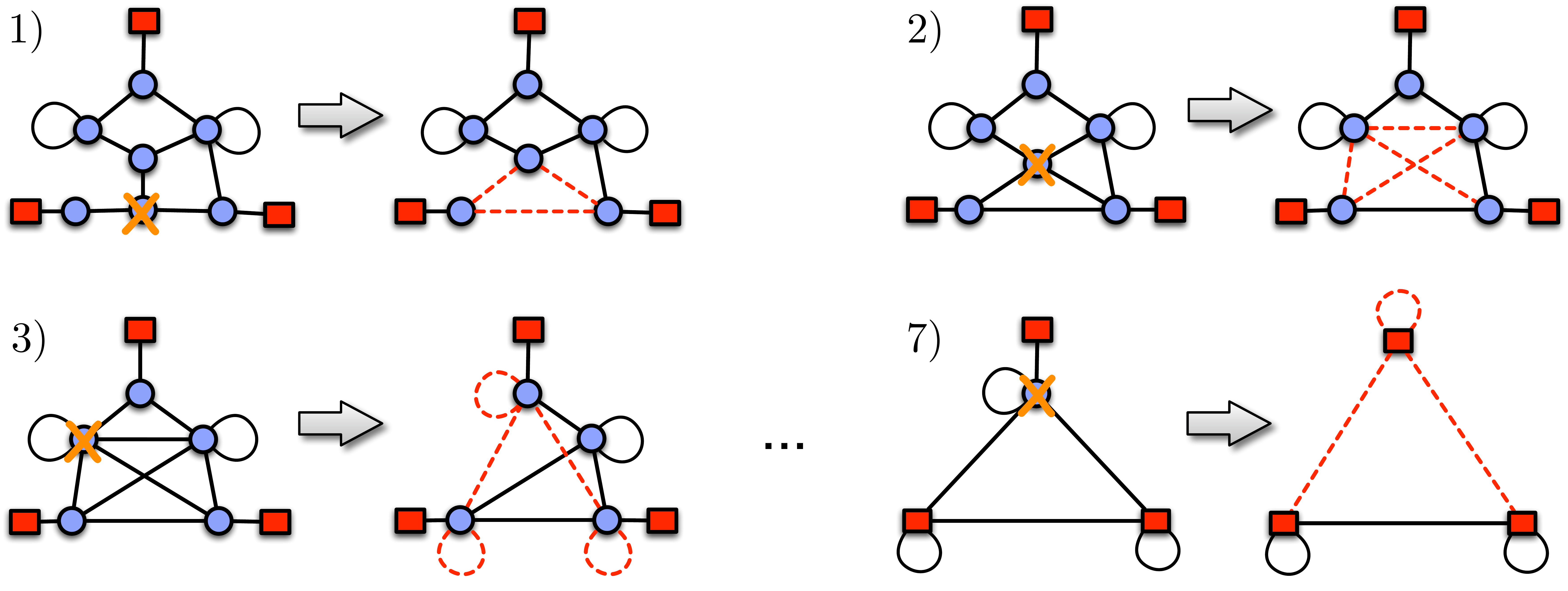}
	\caption{Sparsity pattern (or unweighted graph/topolgical evolution) corresponding to the iterative Kron reduction \eqref{eq: iterative Kron reduction in vector form} of a simple graph with 3 boundary nodes \generator\!\! and 7 interior nodes \loadbus\!\!. The dashed red lines indicate the newly added edges in a reduction step. Steps 4) - 6) are omitted.}
	\label{Fig: Iterative Kron reduction}
	}
\end{figure}

If the sequence \eqref{eq: iterative Kron reduction in vector form} is well-defined, then each $Q^{\ell}$ is a loopy Laplacian matrix inducing a graph by Lemma \ref{Lemma: Structural Properties of Kron Reduction}. Before going further into the details of iterative Kron reduction, we illustrate the unweighted graph corresponding to $Q^{\ell}$ (the sparsity pattern of the corresponding adjacency matrix) in Figure \ref{Fig: Iterative Kron reduction}.
The following observations can be made: (1) The connectivity is maintained. (2) At the $\ell$th reduction step a new edge between two nodes appears if and only if both were connected to $k_{\ell}$ before the reduction, and (3) all other edges persist. (4) Likewise, a new self-loop appears at a node $i \neq k_{\ell}$ if and only if $i$ was connected to $k_{\ell}$ and $k_{\ell}$ featured a self-loop before the reduction. The next subsection turns these observations into rigorous theorems. 

In components, $Q^{\ell}$ is defined according to the celebrated Kron reduction formula
\begin{equation}
	Q_{ij}^{\ell} 
	=
	Q_{ij}^{\ell-1} - \frac{Q^{\ell-1}_{ik_{\ell}} \, Q^{\ell-1}_{jk_{\ell}}}{Q^{\ell-1}_{k_{\ell}k_{\ell}}}
	\label{eq: iterative Kron reduction in components}
	\,,
	\quad
	i,j \in \until {n-\ell}\,,
\end{equation}
which nicely illustrates the step-wise Gaussian elimination. For a well-defined sequence $\{Q^{\ell}\}_{\ell=1}^{n-|\alpha|}$, we let $A^{\ell}$ and $L^{\ell}$ be the corresponding adjacency and loop-less Laplacian matrix of the $\ell$th reduction step.  
The following lemma states some important properties of iterative Kron reduction. In particular, the iterative Kron reduction is well-posed and ultimately results in  the Kron-reduced matrix.

\begin{lemma}[\bf Properties of iterative Kron reduction]
\label{Lemma: Properties of iterative Kron reduction}
Consider the matrix sequence $\{Q^{\ell}\}_{\ell =1}^{n-|\alpha|}$ defined
via iterative Kron reduction in equation~\eqref{eq: iterative Kron
  reduction in vector form}. The following statements hold:
\begin{enumerate}
	
	\item[1)] {\bf Well-posedness:}  Each matrix $Q^{\ell}$, $\ell \in \until{n-|\alpha|}$, is well defined, and the classes of loopy, strictly loopy, and loop-less Laplacian matrices are closed throughout the iterative Kron reduction.
	
	\item[2)] {\bf Quotient property:} The Kron-reduced matrix $\subscr{Q}{red}= Q/Q(\alpha,\alpha)$ can be obtained by iterative reduction of all interior nodes $k_{\ell} \in \mc I_{n} \setminus \alpha$, that is, $\subscr{Q}{red} \equiv Q^{n-|\alpha|}$. Equivalently, the following diagram commutes:
	\begin{center}
		\includegraphics[scale=0.31]{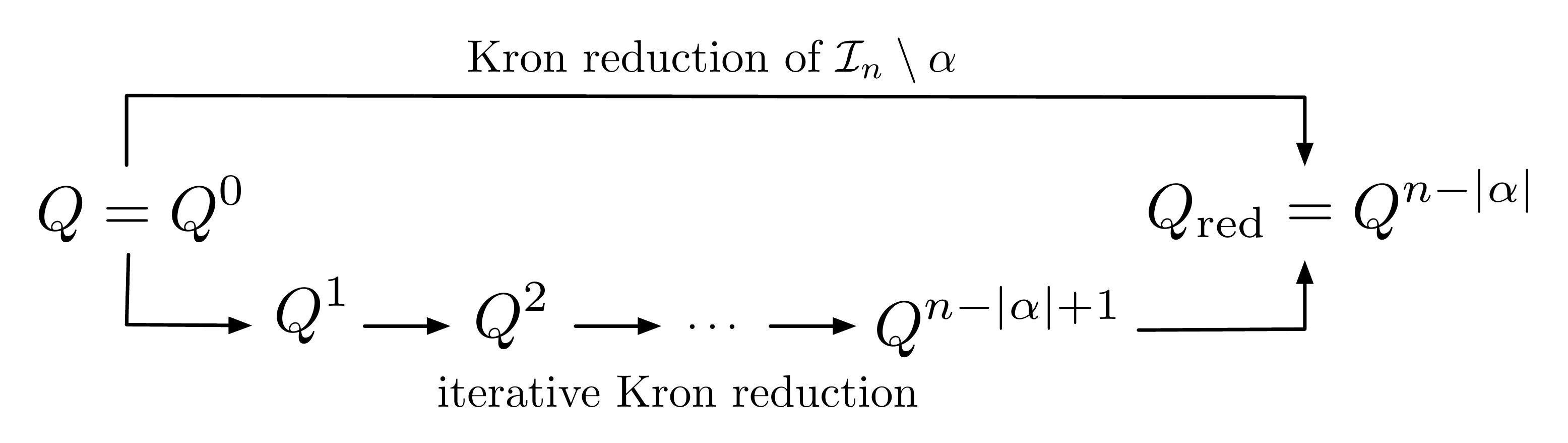}
	\end{center}	

	\item[3)] {\bf Diagonal dominance:} The $i$th row sum of $Q^{\ell}$, $i \in \until{n-\ell}$, is given\,\,by
\begin{equation}
	\sum_{j=1}^{n-\ell} Q^{\ell}_{ij}
	=
	A_{ii}^{\ell}
	=
	\left\{\begin{array}{lr}
	A_{ii}^{\ell-1} \,,
	&\mbox{ if } A^{\ell-1}_{k_{\ell} k_{\ell}} = 0 \,, \\
	A_{ii}^{\ell-1} + A^{\ell-1}_{ik_{\ell}} \Bigl(1 -  \frac{L^{\ell-1}_{k_{\ell}k_{\ell}}}{L^{\ell-1}_{k_{\ell}k_{\ell}} + A^{\ell-1}_{k_{\ell} k_{\ell}}} \Bigr) \,,
	&\mbox{ if } A^{\ell-1}_{k_{\ell} k_{\ell}} > 0 \,.
	\end{array}\right.
	\label{eq: diagonal dominance after one step reduction}
\end{equation}
\end{enumerate}
\end{lemma}

Property 3) of Lemma \ref{Lemma: Properties of iterative Kron reduction} confirms the closure of strictly loopy and loop-less Laplacians and further reveals that the weights of the self-loops are non-decreasing.

{\em Proof of Lemma \ref{Lemma: Properties of iterative Kron reduction}}.
Statement 2) is simply the {\it Quotient Formula} \cite[Theorem 1.4]{FZ:05} stating that Schur complements (or Gaussian elimination for that mater) can be taken iteratively or in a single step. Furthermore, the Quotient Formula states that all intermediate Schur complements $Q^{\ell}$ exist for $\ell \in \until{n-|\alpha|}$. This fact together with the closure properties in Lemma \ref{Lemma: Structural Properties of Kron Reduction} proves statement 1). 
For notational simplicity and without loss of generality, we prove statement 3) for $\ell =1$ and $k_{1} = n$. Note that $A^{0}=A$, $L^{0} = L$, $Q^{0} = Q$, and consider the $i$th row sum of $Q^{1}$ given by
\begin{align*}
	\sum_{j=1}^{n-1} Q^{1}_{ij}
	&=
	\sum_{j=1}^{n-1} \left(Q_{ij} - \frac{Q_{in}\,Q_{jn}}{Q_{nn}} \right)
	=
	\sum_{j=1}^{n-1} \left(Q_{ij} - \frac{A_{in}\,A_{jn}}{L_{nn} + A_{nn}} \right)
	\\
	&=
	\sum_{j=1}^{n-1} \left(Q_{ij} \right) - \frac{A_{in}}{L_{nn} + A_{nn}}  \sum_{j=1}^{n-1} \left( A_{jn} \right)
	=
	A_{ii} + A_{in} - \frac{A_{in}}{L_{nn} + A_{nn}} L_{nn}
	\,,
\end{align*}
where\,we\,used equality \eqref{eq: iterative Kron reduction in components}, the identities $Q \!=\! L+\diag(\{A_{ii}\}_{i=1}^{n})$, $\sum_{j=1}^{n-1} Q_{ij} \!=\!  A_{ii} + A_{in}$, and $\sum_{j=1}^{n-1} A_{jn} = L_{nn}$. Since $A_{nn} \geq 0$ (due to property 1) nonnegative row sums follow also in the general case), we are left with the two cases stated in \eqref{eq: diagonal dominance after one step reduction}.
\qquad\endproof

\subsection{Topological, Spectral, and Algebraic Properties of Kron Reduction}

In this section we finally begin our characterization of the properties of
Kron reduction.  We start by discussing how the graph topology of $G$, or
equivalently the sparsity pattern of $Q$, changes after Kron reduction of
the interior nodes.

\begin{theorem}[\bf Topological Properties of Kron Reduction]
\label{Theorem: Topological Properties}
Let $G$, $\subscr{G}{red}$, and $\widehat G$ be the undirected weighted  graphs associated to $Q$, $\subscr{Q}{red} = Q/Q(\alpha,\alpha)$, and the augmented loopy Laplacian $\widehat Q$, respectively. The following statements hold:
\begin{enumerate}

\item[1)] {\bf Edges:} Two nodes $i,j \in \alpha$ are connected by an edge
  in $\subscr{G}{red}$ if and only if there is a path from $i$ to $j$ in
  $G$ whose nodes all belong to $\{i,j\} \cup ( \mc I_{n} \setminus
  \alpha)$.

\item[2)] {\bf Self-loops:} A node $i \in \alpha$ features a self-loop in
  $\subscr{G}{red}$ if and only if there is a path from $i$ to the grounded
  node $n+1$ in $\widehat G$ whose nodes all belong to $\{i\} \cup ( \mc
  I_{n} \setminus \alpha)$. Equivalently, a node $i \in \alpha$ features a
  self-loop in $\subscr{G}{red}$ if and only if $i$ features a self-loop in
  $G$ or there is a path from $i$ to a loopy interior node $j \in \mc I_{n}
  \setminus \alpha$ whose nodes all belong to $\{i,j\} \cup (\mc I_{n}
  \setminus \alpha)$.
	

\item[3)] {\bf Reduction of connected components:} If the interior nodes
  $\beta \subseteq \mc I_{n} \setminus \alpha$ form a connected subgraph of
  $G$, then the boundary nodes $\bar\alpha \subseteq \alpha$ adjacent to
  $\beta$ in $G$ form a clique in $\subscr{G}{red}$.  Moreover, if one node
  in $\beta$ features a self-loop in $G$, then all boundary nodes adjacent
  to $\beta$ in $G$ feature self-loops in $\subscr{G}{red}$.
\end{enumerate}
\end{theorem}

The topological evolution of the graph corresponding to the iterative Kron
reduction \eqref{eq: iterative Kron reduction in components} -- the
sparsity pattern of $A^{\ell}$ -- is illustrated in Figure \ref{Fig:
  Iterative Kron reduction}. Statement 1) of Theorem \ref{Theorem:
  Topological Properties} can be observed in each reduction step,
statement 2) is visible in the third step, and statement 3) is visible in the final
step of the reduction in Figure \ref{Fig: Iterative Kron reduction}.

In the loop-less case, Theorem \ref{Theorem: Topological Properties} is also proved in \cite[Theorem 4.20 and Theorem 4.23]{BA:02} and in \cite[Theorem 14.2]{MF:86}. In the strictly loopy case, Theorem \ref{Theorem: Topological Properties} can be obtained from \cite[Lemma 3.1]{JJM-MN-HS-MJT:95}. Given our prior results on iterative Kron reduction and the augmented Laplacian matrix, the proof of Theorem \ref{Theorem: Topological Properties} is rather straightforward.

{\em Proof of Theorem \ref{Theorem: Topological Properties}}. 
To prove Theorem \ref{Theorem: Topological Properties} we initially focus on the loop-less case and the reduction of a single interior node $k$ via the one-step iterative Kron reduction \eqref{eq: iterative Kron reduction in components}. Due to the closure of loop-less Laplacian matrices under iterative Kron reduction, c.f. Lemma \ref{Lemma: Properties of iterative Kron reduction}, we can restrict the discussion to the non-positive off-diagonal elements of $Q^{1} \triangleq Q/Q_{kk}$ inducing the mutual edges in the graph. Any non-zero and thus strictly negative element $Q_{ij}$ is rendered to a strictly negative element $Q_{ij}^{1}$ since the first term on the right-hand side of equation \eqref{eq: iterative Kron reduction in components} is strictly negative and the second term is non-positive. Therefore, all edges in the graph induced by $Q_{ij}$ persist in the graph induced by $Q_{ij}^{1}$. According to the iterative Kron reduction formula \eqref{eq: iterative Kron reduction in components}, a zero element $Q_{ij}=0$ is converted into a strictly negative element $Q_{ij}^{1}<0$ if and only if both nodes $i$ and $j$ are adjacent to $k$. Consequently, a reduction of node $k$ leads to a complete graph among all nodes that were adjacent to $k$. In summary, statement 1) holds in the loop-less case and for a single reduction step.

Recall from Lemma \ref{Lemma: Properties of iterative Kron reduction} that the one-step reduction of all interior nodes is equivalent to iterative reduction of each interior node. Hence, the arguments of the previous paragraph can be applied iteratively, which proves statement 1) in the loop-less case.

The results in statement 2), the strictly loopy case, follow simply by applying the same arguments to the augmented Laplacian $\widehat Q$ defined in \eqref{eq: augmented Laplacian}. Alternatively, an element-wise analysis of $A_{ii}^{1}$ together with equation \eqref{eq: diagonal dominance after one step reduction} leads to the same conclusion. 

Finally, statement 3) follows directly by applying statement 1) and 2) to the connected component $\beta$. This completes the proof of Theorem \ref{Theorem: Topological Properties}.
\qquad\endproof 

By Theorem \ref{Theorem: Topological Properties}, the topological connectivity among the boundary nodes becomes only denser under Kron reduction. Hence, the algebraic connectivity $\lambda_{2}(L)$ -- as a spectral connectivity measure -- should increase accordingly. Indeed, for the graph in Figure \ref{Fig: Iterative Kron reduction} (with initially unit weights), we have $\lambda_{2}(L) = 0.30 \leq \lambda_{2}(\subscr{L}{red}) = 0.45$. Physical intuition also suggests that loads in an electrical circuit weaken the influence of nodes on another. Thus, self-loops should weaken the algebraic connectivity $\lambda_{2}(\subscr{L}{red})$ in the reduced network. In the following we confirm these intuitions and present some facts on the spectrum of the Kron-reduced matrix.

\begin{theorem}[\bf Spectral Properties of Kron Reduction]
\label{Theorem: Spectral Properties}
The following statements hold for the spectrum of the Kron-reduced matrix $\subscr{Q}{red} = Q/Q(\alpha,\alpha)$:
\begin{enumerate}
	
	\item[1)] {\bf Spectral interlacing:} For any $r \in \mc I_{|\alpha|}$ it holds that
\begin{equation}
	\lambda_{r}(Q) 
	\leq
	\lambda_{r}(\subscr{Q}{red})
	\leq 
	\lambda_{r}(Q[\alpha,\alpha])
	\leq
	\lambda_{r + n - |\alpha|}(Q)
	\label{eq: bound on Laplacian eigenvalues and their Schur complements}
	\,.
\end{equation}
	
	\item[2)] {\bf Effect of self-loops:} Specifically, for any $r \in \mc I_{|\alpha|}$ it holds that
\begin{align}
	\lambda_{r}(\subscr{L}{red}) + \max_{i \in \alpha} \{ \subscr{A}{red}[i,i] \}
	&\geq 
	\lambda_{r}(L) + \min_{i \in \mc I_{n}} \{ A_{ii} \}
	\label{eq: bound on algebraic connectivity in the loopy case - geq}	
	\,,
	\\
	\lambda_{r}(\subscr{L}{red}) + \min_{i \in \mc I_{|\alpha|}} \{ \subscr{A}{red}[i,i] \}
	&\leq 
	\lambda_{r+n-|\alpha|}(L) + \max_{i \in \mc I_{n}} \{ A_{ii} \}
	\label{eq: bound on algebraic connectivity in the loopy case - leq}
	\,.	
\end{align}
\end{enumerate}
\end{theorem}

In the loop-less case, inequality \eqref{eq: bound on Laplacian eigenvalues and their Schur complements} implies non-decreasing algebraic connectivity $ \lambda_{2}(L) \leq \lambda_{2}(\subscr{L}{red})$. For instance, the graph in Figure \ref{Fig: Example for Kron reduction} with zero-valued self-loops satisfies $\lambda_{2}(L) = 0.39 \leq \lambda_{2}(\subscr{L}{red}) = 0.69$. However, in the strictly loopy case inequalities \eqref{eq: bound on algebraic connectivity in the loopy case - geq}-\eqref{eq: bound on algebraic connectivity in the loopy case - geq} imply that self-loops weaken the algebraic connectivity tremendously. The same graph with unit-valued self-loops satisfies $\lambda_{2}(L) = 0.39 \geq \lambda_{2}(\subscr{L}{red}) = 0.29$. 

To prove Theorem \ref{Theorem: Spectral Properties}, we recall the {\it interlacing property} for Schur complements.

\begin{lemma}\normalfont{({\bf Spectral Interlacing Property}, \cite[Theorem 3.1]{YF:02}).}
\label{Lemma: Interlacing property}
Let $A$ be a Hermitian positive semidefinite matrix of order $n$ and let $\beta$ be a non-empty proper subset of $\mc I_{n}$ such that $A[\beta,\beta]$ is nonsingular. Then for any $r \in \mc I_{n - |\beta|}$, it holds that 
$
	\lambda_{r}(A) 
	\leq
	\lambda_{r}(A/A[\beta,\beta])
	\leq 
	\lambda_{r}(A(\beta,\beta))
	\leq
	\lambda_{r + |\beta|}(A)
$.
\end{lemma}

{\em Proof of Theorem \ref{Theorem: Spectral Properties}}.
Since $Q$ is a loopy Laplacian matrix and hence positive semidefinite, Lemma \ref{Lemma: Interlacing property} can be applied with $\beta =  \mc I_{n} \setminus \alpha$ and results in the bound \eqref{eq: bound on Laplacian eigenvalues and their Schur complements} in statement 1). To prove statement 2), we recall Weyl's inequality \cite[Theorem 4.3.1]{RAH-CRJ:85} for the spectrum of the sum of two Hermitian matrices $A,B$ of dimension $n$:
\begin{equation}
	\lambda_{k}(A) + \lambda_{1}(B) 
	\leq 
	\lambda_{k}(A+B)
	\leq
	\lambda_{k}(A) + \lambda_{n}(B)
	\,,
	\quad
	k \in \mc I_{n}
	\label{eq: Weyl's inequality}
	\,.
\end{equation}
In order to continue, let $r \in \mc I_{|\alpha|}$ (or equivalently, $r \in \alpha$) and consider the following set of spectral equalities and inequalities resulting in the spectral bound \eqref{eq: bound on algebraic connectivity in the loopy case - geq}:
\begin{align*}
	\lambda_{r}(\subscr{L}{red})
	&=
	\lambda_{r}(\subscr{Q}{red} - \diag(\{\subscr{A}{red}[i,i]\}_{i \in \alpha}))
	\geq
	\lambda_{r}(\subscr{Q}{red}) - \max\nolimits_{i \in \alpha} \{ \subscr{A}{red}[i,i] \}
	\\
	&\geq
	\lambda_{r}(Q) - \max\nolimits_{i \in \alpha} \{ \subscr{A}{red}[i,i] \}
	=
	\lambda_{r}(L + \diag(\{A_{ii}\}_{i=1}^{n})) - \max\nolimits_{i \in \alpha} \{ \subscr{A}{red}[i,i] \}
	\\
	&\geq
	\lambda_{r}(L) + \min\nolimits_{i \in \mc I_{n}} \{ A_{ii} \} - \max\nolimits_{i \in \alpha} \{ \subscr{A}{red}[i,i] \}
	\,,
\end{align*}
where we subsequently made use of the identity $\subscr{L}{red} = \subscr{Q}{red} - \diag(\{\subscr{A}{red}[i,i]\}_{i \in \alpha})$, Weyl's inequality \eqref{eq: Weyl's inequality}, the fact $\lambda_{1}(-\diag(\{\subscr{A}{red}[i,i]\}_{i \in \alpha})) = - \max_{i \in \alpha} \{ \subscr{A}{red}[i,i] \}$, Lemma \ref{Lemma: Interlacing property}, the identity $Q = L + \diag(\{A_{ii}\}_{i=1}^{n})$, and again Weyl's inequality \eqref{eq: Weyl's inequality} with $\lambda_{1}(\diag(\{A_{ii}\}_{i=1}^{n})) = \min_{i \in \mc I_{n}} \{ A_{ii} \}$. The bound \eqref{eq: bound on algebraic connectivity in the loopy case - leq} follows analogously. 
\qquad\endproof 

The following theorem summarizes some algebraic properties of the Kron reduction process. In particular, it quantifies the topological properties stated in Theorem \ref{Theorem: Topological Properties}, it quantifies the reduced self-loops occurring in Theorem \ref{Theorem: Spectral Properties}, and it states that the edge and self-loop weights among the boundary nodes are non-decreasing, as seen in Figure \ref{Fig: Example for Kron reduction}. Furthermore, Theorem \ref{Theorem: Algebraic Properties} confirms the intuition that the class of undirected connected graphs is closed under Kron reduction due to closure of irreducibility, and it reveals some more subtle properties concerning the effect of self-loops.

\begin{theorem}
[\bf Algebraic Properties of Kron Reduction]
\label{Theorem: Algebraic Properties}
Consider the Kron-reduced matrix $\subscr{Q}{red}$ and the accompanying matrices $\subscr{Q}{ac} = -Q[\alpha,\alpha)Q(\alpha,\alpha)^{-1}$ and $\subscr{L}{ac} = -L[\alpha,\alpha)L(\alpha,\alpha)^{-1}$ defined in Lemma \ref{Lemma: Structural Properties of Kron Reduction}. 
The following statements hold:
\begin{enumerate}
	
	\item[1)] {\bf Closure of irreducibility:} If $Q$ is irreducible, then $\subscr{Q}{red}$ is irreducible. 
	
	\item[2)] {\bf Monotonicity of elements:} For all $i,j \in \alpha$
          it holds that $\subscr{Q}{red}[i,j] \leq
          Q[\alpha,\alpha][i,j]$. Equivalently, $\subscr{A}{red}[i,j] \geq
          A[\alpha,\alpha][i,j]$ for all $i,j \in \alpha$.
	
	\item[3)] {\bf Effect of self-loops I:} Define $\Delta_{i}
          \triangleq A_{ii} \geq 0$, for $i \in \mc I_{n}$, so that loopy
          and loop-less Laplacians $Q$ and $L$ are related by $Q = L +
          \diag(\{\Delta_{i}\}_{i \in \mc I_{n}})$. Then the Kron-reduced
          matrix takes the form
	\begin{equation}
		\subscr{Q}{red} 
		=
		L/L(\alpha,\alpha) 
		+ \diag(\{\Delta_{i}\}_{i \in \alpha}) + S
		\label{eq: complete reduction of all loopy nodes}
		\,,
	\end{equation}	
	where $S = \subscr{L}{ac}
		( I_{n-|\alpha|} + \diag(\{\Delta_{i}\}_{i \in \mc I_{n}\setminus\alpha}) L(\alpha,\alpha)^{-1} )^{-1}
		 \diag(\{\Delta_{i}\}_{i \in \mc I_{n}\setminus \alpha}) \subscr{L}{ac}^{T}
	$ is a symmetric nonnegative $|\alpha| \times |\alpha|$ matrix. Furthermore, the reduced self-loops satisfy $\subscr{A}{red}[i,i] = \Delta_{i} + \sum_{j=1}^{n-|\alpha|} \subscr{Q}{ac}[i,j] \Delta_{|\alpha|+j}$ for $i \in \alpha$. 
	
	\item[4)] {\bf Effect of self-loops II:} If the subgraph among the interior nodes $\mc I_{n} \setminus \alpha$ is connected, each boundary node $\alpha$ is connected to at least one interior node in $\mc I_{n} \setminus \alpha$, and at least one of the interior nodes has a positively weighted self-loop, then $S$ and $\subscr{Q}{ac}$ are both positive matrices.
	
\end{enumerate}
\end{theorem}

Statements 1) and 2) are not surprising given our knowledge from Theorems
\ref{Theorem: Topological Properties} and \ref{Theorem: Spectral
  Properties}. Statement 3) reveals an interesting fact that can be nicely
illustrated by considering the reduction of a single interior node $k$ with
a self-loop $\Delta_{k}\geq0$. In this case, the matrix $S$ in identity
\eqref{eq: complete reduction of all loopy nodes} specializes to the
symmetric and nonnegative matrix $S = c_{k} \cdot L(k,k] L[k,k) \in
\mbb R^{(n-1) \times (n-1)}$, where $c_{k} = \Delta_{k}/(L_{kk}(L_{kk} +
\Delta_{k})) \geq 0$.
Hence, the reduction of node $k$ decreases the mutual coupling $\{i,j\}$ in
$Q/Q_{kk}$ by the amount $c_{k} \cdot A_{ik} \, A_{jk} > 0$ and increases
each self-loop $i$ in $Q/Q_{kk}$ by the corresponding amount $c_{k} \cdot
A_{ik} \, A_{ik} > 0$. This argument can also be applied iteratively. The
complete reduction of all nodes in statement 4) implies, under quite
general topological conditions, that a single positive self-loop in the
interior network will affect the entire reduced network by decreasing all
mutual connections and increasing all self-loops\,\,weights.

For the proof of Theorem \ref{Theorem: Algebraic Properties}, we recall the following identities for the inverse of the sum of two matrices:

\begin{lemma}\normalfont{({\bf Sherman-Morrison Formula}, \cite{HVH-SRS:81}).}
\label{Lemma: Identities for the inverse of a sum of matrices}
Let $A,B \in \mbb C^{n \times n}$ and let $A$ be non-singular. If $A+B$ is nonsingular, then 
\begin{align}
	( A + B )^{-1}
	&= A^{-1} - A^{-1} (I + BA^{-1})^{-1} BA^{-1} 
	\label{eq: identity for inv(A+B) - 1} 
	\,.
\end{align}
If additionally $B = \Delta \cdot u v^{T}$ for $\Delta \in \mbb R$ and $u,v \in \mbb R^{n}$, then
\begin{align}
	( A + \Delta \cdot u v^{T} )^{-1}
	&=
	A^{-1} - \frac{\Delta}{1 + \Delta \cdot v^{T} A^{-1} u} A^{-1} u v^{T} A^{-1}
	\,.
	\label{eq: identity for inv(A+B) - 4}
\end{align}
\end{lemma}

{\em Proof of Theorem \ref{Theorem: Algebraic Properties}}.
In the loop-less case, the spectral inequality \eqref{eq: bound on Laplacian eigenvalues and their Schur complements} in Theorem \ref{Theorem: Spectral Properties} implies non-decreasing algebraic connectivity $\lambda_{2}(\subscr{L}{red}) \geq \lambda_{2}(L)>0$ and thus irreducibility of $\subscr{L}{red}$. In the strictly loopy case, note that the Kron-reduced graph features the same edges (excluding self-loops) as in the loop-less case, by Theorem \ref{Theorem: Topological Properties}. Thus, connectivity and irreducibility follow, which proves statement 1).

The element-wise bound $\subscr{Q}{red}[i,j] \leq Q[\alpha,\alpha][i,j]$ follows directly from \cite[Lemma 1]{DEC:66}, where this bound is stated for the reduction of one node. By Lemma \ref{Lemma: Properties of iterative Kron reduction}, a one-step reduction is equivalent to iterative one-dimensional reductions. Hence, \cite[Lemma 1]{DEC:66} can be applied iteratively and yields $\subscr{Q}{red}[i,j] \leq Q[\alpha,\alpha][i,j]$. For the off-diagonal elements $i \neq j$, this bound is readily converted to $\subscr{A}{red}[i,j] \geq A[\alpha,\alpha][i,j]$. The same bound follows for the diagonal elements since diagonal dominance is non-decreasing under Kron reduction, see Lemma \ref{Lemma: Properties of iterative Kron reduction}. This completes the proof of statement 2). 

Identity \eqref{eq: complete reduction of all loopy nodes} in statement 3) follows by expanding the Kron-reduced matrix as
\begin{align*}
	\subscr{Q}{red}
	&=
	Q/Q(\alpha,\alpha)
	=
	\bigl( L + \diag(\{\Delta_{i}\}_{i \in \mc I_{n}})\bigr) / \bigl( L(\alpha,\alpha) + \diag( \{\Delta_{i}\}_{i \in \mc I_{n}\setminus \alpha}) \bigr)
	\\
	&=
	\diag(\{\Delta_{i}\}_{i \in \alpha}) + L[\alpha,\alpha]
	- L[\alpha,\alpha) \bigl(L(\alpha,\alpha) + \diag(\{\Delta_{i}\}_{i \in \mc I_{n}\setminus\alpha}) \bigr)^{-1}  L(\alpha,\alpha]
	\\
	&=
	L/L(\alpha,\alpha) + \diag(\{\Delta_{i}\}_{i \in \alpha}) + S
	\,,
\end{align*}	
where $S$ is defined statement 3), and we applied the matrix identity \eqref{eq: identity for inv(A+B) - 1} in the fourth inequality with $A = L(\alpha,\alpha)$ and $B = \diag(\{\Delta_{i}\}_{i \in \mc I_{n}\setminus\alpha})$. 
By Lemma \ref{Lemma: Properties of iterative Kron reduction}, the one-step Schur complement $Q/Q(\alpha,\alpha)$ is equivalent to iterative one-dimensional reduction of all interior nodes $\mc I_{n} \setminus \alpha$, and the matrix $Q^{\ell} = Q^{\ell-1}/ Q^{\ell-1}_{k_{\ell}k_{\ell}}$ at the $\ell$th reduction step is again a loopy Laplacian. If we abbreviate the self-loops at the $\ell$th reduction step by $\Delta_{i}^{\ell} \triangleq A_{ii}^{\ell}$, then  $Q^{\ell}$ can be reformulated according to identity \eqref{eq: complete reduction of all loopy nodes}\,\,as
\begin{equation}
	Q^{\ell}
	=
	Q^{\ell-1}/ Q^{\ell-1}_{k_{\ell}k_{\ell}} 
	=
	L^{\ell-1}/ L^{\ell-1}_{k_{\ell}k_{\ell}} + \diag(\{\Delta^{\ell}_{i}\}_{i=1}^{n-\ell}) + S^{\ell}
	\label{eq: iterative reduction of loopy nodes}
	\,,
\end{equation}
where $S^{\ell} = c_{\ell} \cdot L^{\ell-1}(k_{\ell},k_{\ell}] L^{\ell-1}[k_{\ell},k_{\ell}) = c_{\ell} \cdot A^{\ell-1}(k_{\ell},k_{\ell}] A^{\ell-1}[k_{\ell},k_{\ell})$ is a symmetric and nonnegative matrix and $c_{\ell} = \Delta^{\ell}_{k}/(L_{k_{\ell}k_{\ell}}(L_{k_{\ell}k_{\ell}} + \Delta^{\ell}_{k_{\ell}})) \geq 0$. Iterative application of this argument implies that $S$ must also be a symmetric and nonnegative matrix.

To obtain an explicit expression for the reduced self-loops, re-consider the identity \eqref{eq: Equation determining row sums of Q} defining the self-loops of $L$. In the general loopy case identity \eqref{eq: Equation determining row sums of Q} reads as $Q \fvec 1_{n} = \Delta$.
Block-Gaussian elimination of the interior nodes yields
$
\subscr{Q}{red} \fvec 1_{|\alpha|} 
= \Delta[\alpha] + \subscr{Q}{ac} \Delta(\alpha)
$. Hence, the $i$th row sum of $\subscr{Q}{red}$ satisfies
$\subscr{A}{red}[i,i] = \sum_{j=1}^{|\alpha|} \subscr{Q}{red}[i,j] = \Delta_{i} + \sum_{j=1}^{n-|\alpha|} \subscr{Q}{ac}[i,j] \Delta_{|\alpha|+j}$. This completes the proof of statement 3).

Under the assumptions of statement 4), the positivity of $\subscr{Q}{ac}$ follows from Lemma \ref{Lemma: Structural Properties of Kron Reduction}. To prove positivity of $S$, note that iterative reduction of all but one interior node yields one remaining interior node $k_{n-|\alpha|+1} \triangleq h$.  According to equality \eqref{eq: iterative reduction of loopy nodes}, reduction of this last loopy node gives the matrix $S^{h} =  c_{h} \cdot A^{h}(h,h] A^{h}[h,h)$. Under the assumptions of statement 4), Theorem \ref{Theorem: Topological Properties} implies that $h$ features a self-loop and is connected to all boundary nodes. It follows that $c_{h} > 0$ and $A^{h}_{i h} > 0$ for all $i \in \mc I_{|\alpha|+1}$. Therefore, $S_{h}$ is a positive matrix, and the same can be concluded for $S$. 
\qquad\endproof

\subsection{Kron Reduction and Effective Resistance}

The physical intuition behind the Kron reduction and the effective resistance in Remark \ref{Remark: physical intuition about resistance} suggests that the transfer conductances $\subscr{Q}{red}[i,j]$ are related to the corresponding effective conductances $1/R_{ij}$. The following theorem gives the exact relation between the Kron-reduced matrix $\subscr{Q}{red}$, the effective resistance matrix $R$, as well as the augmented Laplacian $\widehat Q$.

\begin{theorem}[\bf Resistive Properties of Kron Reduction]
\label{Theorem: Resistive Properties of Kron Reduction} 
Consider the Kron-reduced matrix $\subscr{Q}{red} = Q/Q(\alpha,\alpha)$, the effective resistance matrix $R$, and the augmented Laplacian matrix $\widehat Q$ defined in \eqref{eq: augmented Laplacian}. The following statements hold:

\begin{enumerate}

	\item[1)] {\bf Invariance of under Kron reduction:} The effective resistance $R_{ij}$ between any two boundary nodes $i,j \in \alpha$ is equal when computed from $Q$ or $\subscr{Q}{red}$:
\begin{equation}
	R_{ij} =  (e_{i} - e_{j})^{T} Q^{\dagger} (e_{i} - e_{j}) \equiv (e_{i} - e_{j})^{T} \subscr{Q}{red}^{\dagger} (e_{i} - e_{j}) 	
	\,,\quad i,j \in \alpha
	\label{eq: invariance of effective resistance}
	\,.
\end{equation}

	\item[2)] {\bf Invariance under augmentation:} If $Q$ is a strictly loopy Laplacian, then the effective resistance $R_{ij}$ between any two nodes $i,j \in \mc I_{n}$ is equal when computed from $Q$ or $\widehat Q$:
\begin{equation}
	R_{ij}
	=
	(e_{i}-e_{j})^{T} Q^{-1} (e_{i}-e_{j})
	\equiv
	(e_{i}-e_{j})^{T}  
	\widehat Q^{\dagger}
	(e_{i}-e_{j})
	\,,\quad i,j \in \mc I_{n}
	\,.
	\label{eq: effect of self-loops on effective resistance matrix}
\end{equation}
\end{enumerate}
In other words, statements 1) and 2) imply that, if $Q$ is a strictly loopy
Laplacian, then the following diagram commutes:
\begin{center}
		\includegraphics[scale=0.31]{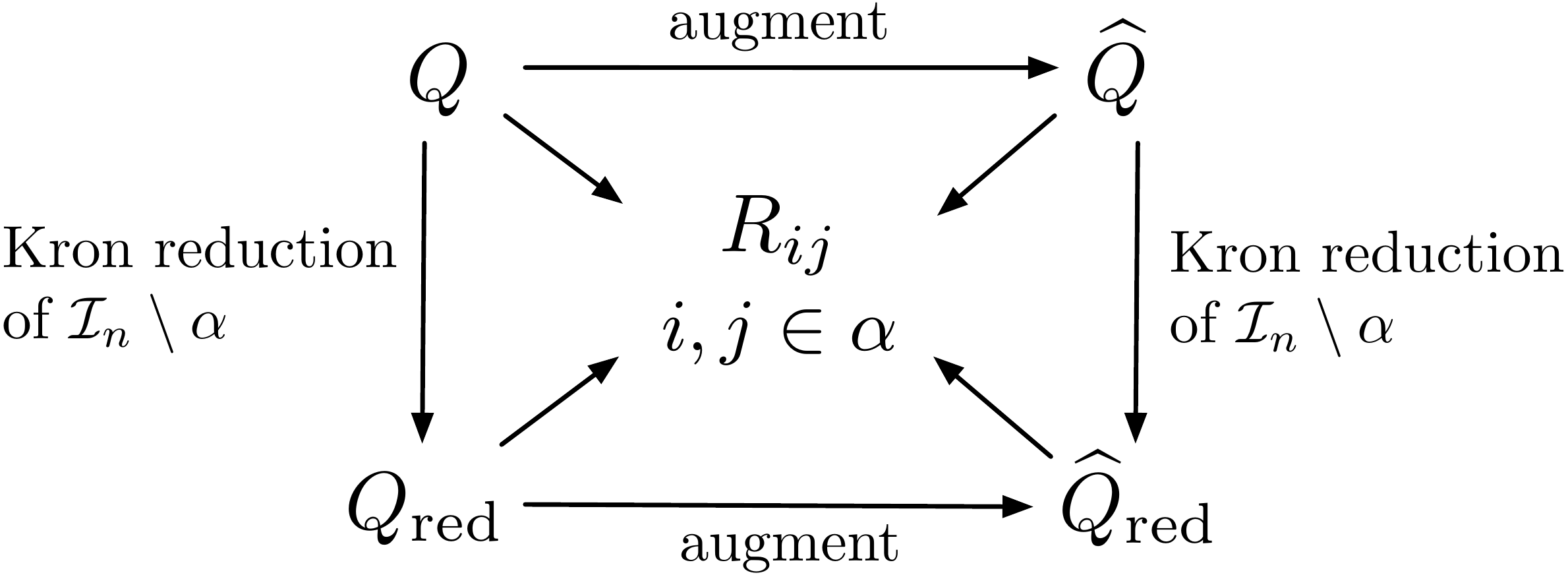}
\end{center}	

\begin{enumerate}

	\item[3)] {\bf Effect of self-loops:} If $Q$ is a strictly loopy Laplacian and $\bar R_{ij} \triangleq (e_{i}-e_{j})^{T} L^{\dagger} (e_{i}-e_{j})$, $i,j \in \mc I_{n}$, is the effective resistance computed from the corresponding loop-less Laplacian $L$, then for all $i,j \in \mc I_{n}$ it holds that $R_{ij} \leq \bar R_{ij}$.
	
\end{enumerate}
\end{theorem}

\begin{figure}[htbp]
	\centering{
	\includegraphics[scale=0.25]{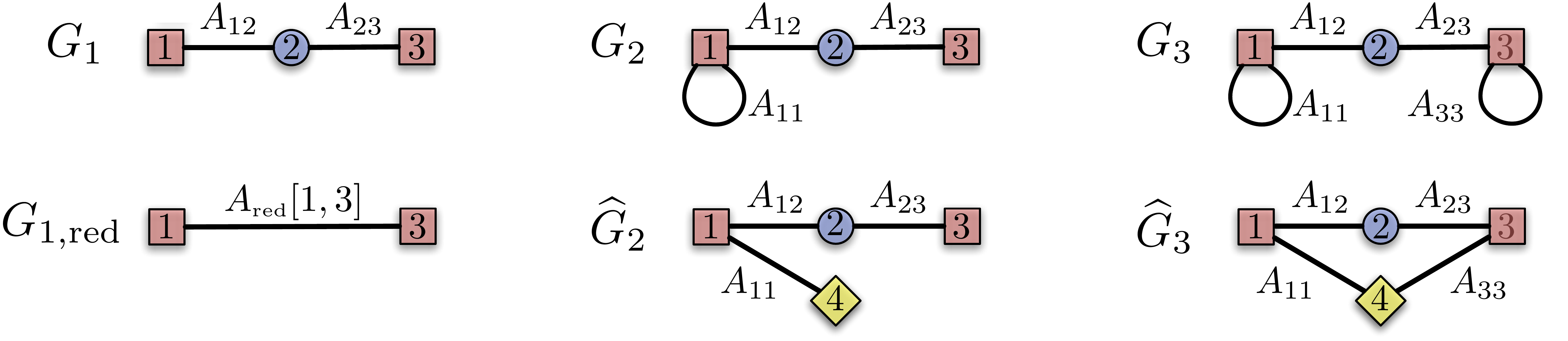}
	\caption{Illustration of Theorem \ref{Theorem: Resistive Properties of Kron Reduction} with six graphs: According to statement 1), the effective resistance $R_{13}$ between the boundary nodes \generator\!\!\! is equal when computed in the graph $G_{1}$ or in the corresponding Kron-reduced graph $G_{1,\textup{red}}$. According to statement 2), the effective resistance $R_{13}$ is equal when computed in the strictly loopy graph $G_{2}$ (resp. $G_{3}$) or in the corresponding augmented loop-less graph $\widehat G_{2}$ (resp. $\widehat G_{3}$) with additional grounded node \terminalbus\!\!. According to statement 3), the effective resistance $R_{13}$ in the strictly loopy graphs $G_{2}$ and $G_{3}$ is not larger than in the loop-less graph $G_{1}$ (with equality for the pair $\{G_{1},G_{2}\}$ and strict inequality for the pair $\{G_{1},G_{3}\}$).}
	\label{Fig: Illustration of the effect of self-loops on the effective resistance}
	}
\end{figure}

Figure \ref{Fig: Illustration of the effect of self-loops on the effective resistance} illustrates the ideas of Theorem \ref{Theorem: Resistive Properties of Kron Reduction} with six simple graphs. 

Identity \eqref{eq: invariance of effective resistance} states that the effective resistances between the boundary nodes are invariant under Kron reduction of the interior nodes. Spoken in terms of circuit theory, the effective resistance between the terminals $\alpha$ can be obtained from either the impedance matrix $Q^{\dagger}$ or the transfer impedance matrix $\subscr{Q}{red}^{\dagger}$. 
Identity \eqref{eq: effect of self-loops on effective resistance matrix} gives a resistive interpretation of the self-loops: the effective resistance among the nodes in a strictly loopy graph $G$ is equivalent to the effective resistance among the corresponding nodes in the augmented loop-less graph $\widehat G$ with additional grounded node. The commutative diagram in Theorem \ref{Theorem: Resistive Properties of Kron Reduction} combines the invariance identity \eqref{eq: invariance of effective resistance} with the self-loop identity \eqref{eq: effect of self-loops on effective resistance matrix}.
According to statement 3), the self-loops do not increase the effective resistance. This is in accordance with the physical interpretation in Remark \ref{Remark: physical intuition about resistance} when viewing the self-loops as shunt conductances. 

To prove Theorem \ref{Theorem: Resistive Properties of Kron Reduction}, we need the following identities relating $R$ and $L$.

\begin{lemma}[\bf Alternative Derivation of the Effective Resistance]
\label{Lemma: Effective resistance obtained differently}
Let $L \in \mbb R^{n \times n}$ be a symmetric irreducible loop-less Laplacian matrix. Then for any $\delta \neq 0$ 
\begin{equation}
	\left(L + (\delta/n) \, \fvec 1_{n \times n} \right)^{-1}
	=
	L^{\dagger} + (1/\delta n) \, \fvec 1_{n \times n}
	\label{eq: matrix identity with pseudo inverse}
	\,.
\end{equation}
Consider the matrix $R$ of effective resistances defined in \eqref{eq: Definition of effective resistance}. For $\delta \neq 0$ it holds that
\begin{equation}
	R_{ij}  
	= (e_{i} - e_{j})^{T} L^{\dagger} (e_{i} - e_{j}) 
	\equiv 
	(e_{i} - e_{j})^{T} ( L + (\delta/n)\, \fvec 1_{n \times n} )^{-1} (e_{i} - e_{j})
	\label{eq: effective resistance from auxiliary matrix}
	\,,\quad i,j \in \mc I_{n}
	\,.
\end{equation}
If further $n \geq 3$, then by taking node $n$ as reference the following identity holds:
\begin{equation}
	R_{ij} 
	= (e_{i} - e_{j})^{T} L^{\dagger} (e_{i} - e_{j}) 
	\equiv 
	(e_{i} - e_{j})^{T} L(n,n)^{-1} (e_{i} - e_{j})
	\label{eq: effective resistance from Dirichlet Laplacian}
	\,,\quad i,j \in \mc I_{n-1}
	\,.
\end{equation}
\end{lemma}

\begin{proof}
Identity \eqref{eq: matrix identity with pseudo inverse} is a simple generalization of \cite[Theorem 5]{IG-WX:04}. Since $\fvec 1_{n \times n}  \fvec 1_{n \times n} = n \cdot \fvec 1_{n \times n}$ and $L L^{\dagger} = L^{\dagger} L = I_{n} - (1/n) \cdot \fvec 1_{n \times n}$ (by definition of $L^{\dagger}$ via the singular value decomposition), identity \eqref{eq: matrix identity with pseudo inverse} can be verified since the product of the left-hand and the right-hand side of \eqref{eq: matrix identity with pseudo inverse} equal the identity matrix. The identity \eqref{eq: effective resistance from auxiliary matrix} follows then by multiplying equation \eqref{eq: matrix identity with pseudo inverse} from the left by $(e_{i} - e_{j})^{T}$ and from the right by $(e_{i} - e_{j})$.
Let $\bar L := L(n,n)$, then it follows from \cite[Appendix B, eq. (17)]{FF-AP-JMR-MS:07} that
$
\bar L^{-1}_{ij} 
=
L^{\dagger}_{ij} - L^{\dagger}_{in} - L^{\dagger}_{jn} - L^{\dagger}_{nn}
$. 
The identity \eqref{eq: effective resistance from Dirichlet Laplacian} can now be directly verified. For an alternative derivation of identity \eqref{eq: effective resistance from Dirichlet Laplacian}, we refer to \cite{AG-SB-AS:08}.
\qquad
\end{proof}

{\em Proof of Theorem \ref{Theorem: Resistive Properties of Kron Reduction}}.
We begin by proving statement 1) in the strictly loopy case when $Q$ is invertible and $Q^{\dagger} = Q^{-1}$. Note that we are interested in the effective resistances only among the nodes $\alpha$, i.e., the $|\alpha| \times |\alpha|$ block of $Q^{-1}$. The Schur complement formula \cite[Theorem 1.2]{FZ:05} gives the $|\alpha| \times |\alpha|$ block of $Q^{-1}$ as $(Q/Q(\alpha,\alpha))^{-1} = \subscr{Q}{red}^{-1}$. Consequently, for $i,j \in \alpha$ the defining equation \eqref{eq: Definition of effective resistance} for the effective resistance $R_{ij}$ is simply rendered to
$ R_{ij}
=
(e_{i} - e_{j})^{T} \subscr{Q}{red}^{-1} (e_{i} - e_{j})
$,
the claimed identity \eqref{eq: invariance of effective resistance}.

In the loop-less case when $Q \equiv L$ is singular, a similar line of arguments holds on the image of $L$. Let $\delta > 0$ and consider the modified and non-singular Laplacian $\tilde L \triangleq L + (\delta / n) \fvec 1_{n \times n}$.  Due to identity \eqref{eq: matrix identity with pseudo inverse} we have that $\tilde L^{-1} = L^{\dagger} + (1/ \delta n) \fvec 1_{n \times n}$. Since $(e_{i} - e_{j})^{T} \fvec 1_{n \times n} (e_{i} - e_{j}) = 0$, we can rewrite identity \eqref{eq: effective resistance from auxiliary matrix} in expanded form as
\begin{equation}
	R_{ij}
	=
	(e_{i} - e_{j})^{T} (L^{\dagger} + (1 /\delta n) \fvec 1_{n \times n})  (e_{i} - e_{j})
	=
	(e_{i} - e_{j})^{T} \tilde L^{-1} (e_{i} - e_{j}) 
	\label{eq: effective resistance from Laplacian with self-loops - 1}
	\,.
\end{equation}
As before, the $|\alpha| \times |\alpha|$ block of $\tilde L^{-1}$ is given by $\tilde L^{-1}$ as $(\tilde L/ \tilde L(\alpha,\alpha))^{-1}$ \cite[Theorem 1.2]{FZ:05}. Consequently, the identity \eqref{eq: effective resistance from Laplacian with self-loops - 1} is for $i,j \in \alpha$ rendered to
\begin{equation}
	R_{ij}
	=
	(e_{i} - e_{j})^{T} (\tilde L/ \tilde L(\alpha,\alpha))^{-1} (e_{i} - e_{j})
	\,.
	\label{eq: effective resistance from Laplacian with self-loops - 2}
\end{equation}
Note that the right-hand side of \eqref{eq: effective resistance from Laplacian with self-loops - 1}, or equivalently \eqref{eq: effective resistance from Laplacian with self-loops - 2}, is independent of $\delta$ since the matrices are evaluated on the subspace orthogonal to $\fvec 1_{n}$, the nullspace of $\tilde L$ as $\delta \downarrow 0$. Thus, on the image of $L$ the limit of the right-hand side of \eqref{eq: effective resistance from Laplacian with self-loops - 2} exists as $\delta \downarrow 0$. By definition, $L^{\dagger}$ acts as regular inverse on the image of $L$, and equation \eqref{eq: effective resistance from Laplacian with self-loops - 2} is rendered to $R_{ij} = (e_{i} - e_{j})^{T} (L/L(\alpha,\alpha))^{\dagger} (e_{i} - e_{j})$. Finally, recall that $\subscr{L}{red} = L/L(\alpha,\alpha)$ which yields the claimed identity \eqref{eq: invariance of effective resistance} in the loop-less case. 

To prove statement 2), note that the strictly loopy Laplacian $Q$ is invertible (due to irreducible diagonal dominance \cite[Corollary~6.2.27]{RAH-CRJ:85}). Hence, the defining equation \eqref{eq: Definition of effective resistance} for the resistance features a regular inverse. The matrix $Q$ can also be seen as the principal $n \times n$ block of the augmented Laplacian $\widehat Q$, that is, $Q = \widehat Q(n+1,n+1)$. The identity \eqref{eq: effect of self-loops on effective resistance matrix} follows then directly from identity \eqref{eq: effective resistance from Dirichlet Laplacian} (with $n$ replaced by $n+1$).

To prove statement 3), we appeal to Rayleigh's celebrated {\em monotonicity law and short/cut principle} \cite{PGD-JLS:84}. Since the Laplacian $L$ induces the same graph as $\widehat Q$ with node $n+1$ cut off, the monotonicity law states that the effective resistance $\bar R_{ij}$ in the graph induced by $L$ is not smaller than the effective resistance $R_{ij}$ in the graph induced by $\widehat Q$. The latter again equals the effective resistance in the graph induced by $Q$ due to identity \eqref{eq: effect of self-loops on effective resistance matrix}. Equivalently, for $i,j \in \mc I_{n}$ it holds that 
$\bar R_{ij} = (e_{i}-e_{j})^{T} L^{\dagger} (e_{i}-e_{j}) \geq (e_{i}-e_{j})^{T} \widehat Q^{\dagger} (e_{i}-e_{j}) = (e_{i}-e_{j})^{T} Q^{-1} (e_{i}-e_{j}) = R_{ij}$.
\qquad\endproof

Theorem \ref{Theorem: Resistive Properties of Kron Reduction} allows to compute the effective resistance matrix $R$ from the transfer impedance matrix $\subscr{Q}{red}^{\dagger}$. We are now interested in a converse result. Iterative\,\,methods constructing $\subscr{Q}{red}^{\dagger}$ from $R$ can be found in \cite{EC-M-JM:94,SN:02,AME:96}. It is also possible to recover the (pseudo) inverse of the loopy Laplacian $Q$, the augmented Laplacian $\widehat Q$, or the corresponding Kron-reduced Laplacians {\em directly} from the effective resistance matrix\,\,$R$.

\begin{lemma}[\bf Impedance and Effective Resistance Matrix]
\label{Lemma: Recovering the Laplacian from the effective resistance}
Let $Q \in \mbb R^{n \times n}$ be a symmetric irreducible loopy Laplacian matrix. Consider the following three cases:
\begin{enumerate}

	\item[1)] {\bf Loop-less case:} Let $R \in \mbb R^{n \times n}$ be the matrix of effective resistances. Then for $i,j \in \mc I_{n}$ the following identity holds:
	\begin{equation}
		Q^{\dagger}[i,j] \equiv L^{\dagger}[i,j] 
		=
		-\frac{1}{2} \Bigl( R_{ij} - \frac{1}{n} \sum_{k=1}^{n} (R_{ik} + R_{jk}) 
		+ \frac{1}{n^{2}} \sum_{k=1}^{n} \sum_{\ell=1}^{n} R_{k\ell} \Bigr)
		\label{eq: Laplacian from effective resistance --- loop-less case}
		\,.
	\end{equation}

	\item[2)] {\bf Strictly loopy case:} Consider the grounded node with index $n+1$, the corresponding augmented Laplacian matrix $\widehat Q \in \mbb R^{(n+1) \times (n+1)}$ defined in \eqref{eq: augmented Laplacian}, and the corresponding matrix of effective resistances $R \in \mbb R^{(n+1) \times (n+1)}$ defined in \eqref{eq: Definition of effective resistance}. Then the following two identities hold:
	\begin{align}
		\widehat Q^{\dagger}[i,j] 
		&=
		-\frac{1}{2} \Bigl( R_{ij} - \frac{1}{n+1} \sum_{k=1}^{n+1} (R_{ik} + R_{jk}) 
		+ \frac{1}{(n+1)^{2}} \sum_{k=1}^{n+1} \sum_{\ell=1}^{n+1} R_{k\ell} \Bigr)
		\,,\nonumber
		\\&\qquad i,j \in \mc I_{n+1}
		\label{eq: Laplacian from effective resistance --- loopy augmented case}
		\,,\\
		Q^{-1}[i,j] 
		&=
		\frac{1}{2}(R[i,n+1] + R[j,n+1] - R[i,j])
		\,,\quad i,j \in \mc I_{n}
		\label{eq: Laplacian from effective resistance --- loopy case}
	\end{align}
	
	\item[3)] {\bf Kron reduced case:} The identities \eqref{eq: Laplacian from effective resistance --- loop-less case}, \eqref{eq: Laplacian from effective resistance --- loopy augmented case}, and \eqref{eq: Laplacian from effective resistance --- loopy case} also hold when $Q^{\dagger}$, $\widehat Q^{\dagger}$, and $Q^{-1}$ on the left-hand sides are replaced by $\subscr{Q^{\dagger}}{red}$, $\subscr{\widehat Q^{\dagger}}{red}$, and $\subscr{Q^{-1}}{red}$, respectively, and $n$ on the right-hand sides is replaced by $|\alpha|$.
	
\end{enumerate}
\end{lemma}

\begin{proof}
Identity \eqref{eq: Laplacian from effective resistance --- loop-less case} is stated in \cite[Theorem 4.8]{EB-AC-AME-JMG:09} for the weighted case and in \cite[Theorem 7]{WXiao-IG:03} for the unweighted case. 
According to statement 2) of Theorem \ref{Theorem: Resistive Properties of Kron Reduction}, the resistance is invariant under augmentation. Hence, identity \eqref{eq: Laplacian from effective resistance --- loop-less case} applied to the augmented Laplacian $\widehat Q$ yields identity \eqref{eq: Laplacian from effective resistance --- loopy augmented case}.
Identity \eqref{eq: Laplacian from effective resistance --- loopy case} follows directly from \cite[Theorem 4.9]{EB-AC-AME-JMG:09} or \cite[Theorem 2.1]{MF:98}. 
According to statement 1) of Theorem \ref{Theorem: Resistive Properties of Kron Reduction}, the effective resistance is invariant under Kron reduction, and thus the effective resistance corresponding to $\subscr{Q}{red}$ is simply $R[\alpha,\alpha]$. Hence, the formulas \eqref{eq: Laplacian from effective resistance --- loop-less case}, \eqref{eq: Laplacian from effective resistance --- loopy augmented case}, and \eqref{eq: Laplacian from effective resistance --- loopy case} can be applied to the Kron-reduced matrix $\subscr{Q}{red}$ as stated in 3).
\qquad
\end{proof}

The formulas stated Theorem \ref{Theorem: Resistive Properties of Kron Reduction} and Lemma \ref{Lemma: Recovering the Laplacian from the effective resistance} allow to compute the effective resistance matrix $R$ in the original non-reduced network from the (pseudo) inverse of the Kron-reduced Laplacian $\subscr{Q}{red}$, and vice versa. In some applications, such as the transient stability problem presented in Subsection \ref{Subsection: transient stability application}, it is desirable to know an {\em explicit} algebraic relationship between $R$ and $\subscr{Q}{red}$ without the (pseudo) inverse.
However, such an explicit relationship between $R$ and $\subscr{Q}{red}$ can be found only when closed-form solutions of $\subscr{Q}{red}^{\dagger}$, $\subscr{Q}{red}^{-1}$ or $\subscr{\widehat Q}{red}^{\dagger}$ are known. These are generally not available.
Generally, it also seems infeasible to relate bounds on $R$ to bounds on $\subscr{Q}{red}$ since element-wise bounding of inverses of interval matrices is known to be NP-hard \cite{GEC:99}. 

Fortunately, closed forms of the inverses of $\subscr{Q}{red}^{\dagger}$, $\subscr{Q}{red}^{-1}$ or $\subscr{\widehat Q}{red}^{\dagger}$ can be derived in an {\it ideal} electric network, where the effective resistances among the boundary nodes are uniform, as well as the effective resistances between the boundary nodes and the ground. In this case, the following theorem states the equivalence of 
such an ideal network and uniform transfer conductances in the Kron-reduced network.

\begin{theorem}[\bf Equivalence of Uniformity in Effective Resistance and Kron Reduction]
\label{Theorem: Effective resistance -- Kron-reduced Laplacian -- uniform resistance}
Consider the Kron-reduced Laplacian $\subscr{Q}{red} = Q/Q(\alpha,\alpha)$ and the corresponding adjacency matrix $\subscr{A}{red}$. Consider the following two cases:

\smallskip
{\bf Loop-less case:} Let $R \in \mbb R^{n \times n}$ be the matrix of effective resistances defined in \eqref{eq: Definition of effective resistance}. Then the following two statements are equivalent:
\begin{enumerate}

	\item[1)] The effective resistances among the boundary nodes $\alpha$ are uniform, i.e., there is $r>0$ such that $R_{ij} = r$ for all distinct $i,j \in \alpha$;
		
	\item[2)] The weighting of the edges in the Kron-reduced network is uniform, i.e., there is $a>0$ such that $\subscr{A}{red}[i,j] = a>0$ for all distinct $i,j \in \alpha$.

\end{enumerate}
If both statements 1) and 2) are true, then it holds that $r = \frac{2}{|\alpha| a}$.

\smallskip

{\bf Strictly loopy case:} Consider additionally the grounded node $n+1$ and the augmented Laplacian matrices $\widehat Q$ and $\subscr{\widehat Q}{red}$ defined in \eqref{eq: augmented Laplacian} and \eqref{eq: augmented Laplacian reduced}, respectively. Let $R \in \mbb R^{(n+1) \times (n+1)}$ be the matrix of effective resistances in the augmented network. Then the following two statements are equivalent:
\begin{enumerate}

	\item[3)] The effective resistances among the boundary nodes $\alpha$ are uniform, i.e., there is $r>0$ such that $R_{ij} = r$ for all distinct $i,j \in \alpha$. The effective resistances between all boundary nodes $\alpha$ and the grounded node $n+1$ are uniform , i.e., there is $g>0$ such that $R[i,n+1] = g$ for all $i \in \alpha$;
	
	\item[4)] The weighting of the edges and the self-loops in the Kron-reduced network is uniform, i.e., there are $a>0$ and $b \geq 0$ such that $\subscr{A}{red}[i,j] = a>0$  and $\subscr{A}{red}[i,i] = b \geq0$ for all distinct $i,j \in \alpha$.

\end{enumerate}
If both statements 3) and 4) are true, then it holds that $r = \frac{2}{|\alpha| a + b}$ and $g= \frac{a+b}{b(a|\alpha|+b)}$.
\end{theorem}

\begin{remark}[Physical networks and graph topologies with uniform resistance]
\normalfont
The assumption of uniform effective resistances in statements 1) and 3) is by no means artificial but rather corresponds to an {\it ideal network}, where all boundary are electrically uniformly distributed in the network, with respect to each other and with respect to the shunt loads. In the applications of electrical impedance tomography and smart grid monitoring, this assumption can be met by choosing the boundary nodes corresponding to sensor locations. In the transient stability problem in power networks the generators corresponding to the boundary nodes are distributed over the power grid ideally in such a way that the loads can be effectively balanced. Otherwise the power grid is vulnerable to voltage collapse and frequency instabilities. Hence, the uniformity assumptions are ideally met in man-made networks.
 
Independently of engineered electrical networks, the following examples show that uniform resistances among a set of boundary nodes $\alpha$ occur for various graph topologies. Here we neglect self-loops and edge weights introducing additional degrees of freedom to construct graphs with uniform resistances. In the trivial case, $|\alpha| = 2$, Theorem \ref{Theorem: Effective resistance -- Kron-reduced Laplacian -- uniform resistance} reduces to \cite[Corollary 4.41]{PETJ-EJPP:09} and the effective resistance among the $\alpha$ nodes is clearly uniform. Second, if the boundary nodes are 1-connected leaves of a highly symmetric graph among the interior nodes, such as a star-shaped tree, a complete graph, or a combination of these two, then the effective resistance among the boundary nodes is uniform. Third, the effective resistance in large-scale small-world networks is uniform among sufficiently distant nodes \cite{GK-MBH-KEB-MJB-BK-DA:06}. Fourth, with increasing number of nodes the effective resistance in random geometric graphs converges to a degree-dependent limit \cite{AR-UL-MH:09}, which is uniform for various geometries and node distributions. Fifth and finally, geometric graphs such as lattices and their fuzzes are special random geometric graphs. According to the previous arguments, the resistance among sufficiently distant lattice nodes becomes uniform in the large limit.
\oprocend
\end{remark}


In order to prove Theorem \ref{Theorem: Effective resistance -- Kron-reduced Laplacian -- uniform resistance}, we need the following Laplacian identities.

\begin{lemma}[\bf Uniform Laplacian Matrices and their Inverses]
\label{Lemma: Laplacian and their inverses}
Let $a>0$ and $b \geq 0$ and consider the loopy $(n \times n)$-dimensional Laplacian matrix
$Q 
\triangleq
a \bigl( n I_{n} - \fvec 1_{n\times n} \bigr) + b I_{n}$ corresponding to a complete graph with $n$ nodes, uniform positive edge weights $a>0$ between any two distinct nodes, and nonnegative and uniform self-loops $b \geq 0$ attached to every node. The following statements hold:
\begin{enumerate}

	\item[1)] For zero self-loops $b=0$, the pseudo inverse of $Q$ is the loop-less Laplacian
\begin{equation*}
	Q^{\dagger}
	=
	\frac{1}{n^{2}a^{2}} \cdot Q
	=
	\frac{1}{n^{2} a} \cdot \bigl( n I_{n} - \fvec 1_{n\times n})
	\,.
\end{equation*}

	\item[2)] For positive self-loops $b>0$, the inverse of $Q$ is the positive matrix
\begin{equation*}
	Q^{-1}
	=
	- \frac{a}{b(an+b)} \bigl( n I_{n} - \fvec 1_{n\times n} \bigr) + \frac{1}{b} I_{n}
	\,.
\end{equation*}

	\item[3)] Consider the augmented Laplacian $\widehat Q$ defined in \eqref{eq: augmented Laplacian} and given as
\begin{equation*}
	\widehat Q
	=
	\left[\begin{array}{c|c}
	a \bigl( n I_{n} - \fvec 1_{n} \bigr) + b I_{n} & - b \fvec 1_{n} 
	\\ \hline
	- b \fvec 1_{n}^{T} & n \cdot b
	\end{array}\right]
	\,.
\end{equation*}
The pseudo inverse of $\widehat Q$ is given by the (augmented) loop-less Laplacian
\begin{equation*}
	\widehat Q^{\dagger}
	=
	\left[\begin{array}{c|c}
	c \bigl( n I_{n} - \fvec 1_{n} \bigr) + d I_{n} & - d \fvec 1_{n} 
	\\ \hline
	- d \fvec 1_{n}^{T} & n \cdot d
	\end{array}\right]
	\,,
\end{equation*}
where $d = 1/(b(n+1)^{2})$ and $c = -d \cdot (a-(n+2)b)/(an+b)$.
\end{enumerate}
\end{lemma}

\begin{proof} 
Lemma \ref{Lemma: Laplacian and their inverses} can be verified by direct computation. 
Statement 2) follows since $Q \cdot Q^{-1} = Q^{-1} \cdot Q = I_{n}$. 
Statements 1)  and 3) follow since $Q$ and $Q^{\dagger}$ (resp. $\widehat Q$ and $\widehat Q^{\dagger}$) satisfy the Penrose equations \cite[Theorem 6.1.1]{DSB:09}. Hence, the loop-less (augmented) Laplacian $Q^{\dagger}$ (resp. $\widehat Q^{\dagger}$) is the unique Moore-Penrose inverse. 
\qquad
\end{proof}

We now have all three ingredients to prove Theorem \ref{Theorem: Effective resistance -- Kron-reduced Laplacian -- uniform resistance}: the invariance formulas \eqref{eq: invariance of effective resistance}-\eqref{eq: effect of self-loops on effective resistance matrix} for the effective resistance stated in Theorem \ref{Theorem: Resistive Properties of Kron Reduction}, the relations between effective resistance and the Kron-reduced impedance matrix in statement 3) of Lemma \ref{Lemma: Recovering the Laplacian from the effective resistance}, and the Laplacian identities in Lemma \ref{Lemma: Laplacian and their inverses}. Given these formulas, the proof of Theorem \ref{Theorem: Effect of rank one perturbation} reduces to mere computation. For the sake of brevity, it will be\,\,omitted.

\subsection{Sensitivity of Kron Reduction to Perturbations}\label{Section: Sensitivity of Kron Reduction to Perturbations}

In the final section of our analysis of Kron reduction we discuss the sensitivity of the Kron-reduced matrix $\subscr{Q}{red}$ to perturbations in the original matrix $Q$.  A number of interesting perturbations can be modeled by adding symmetric matrix $W \in \mbb R^{n \times n}$ and considering the perturbed loopy Laplacian matrix $\tilde Q = Q +W$, where $Q$ is the nominal loopy Laplacian.

The case when $W$ is a diagonal matrix is fully discussed in Theorem \ref{Theorem: Algebraic Properties}. A perturbation of the form when $W[\alpha,\alpha]$ is a non-zero matrix and all other entries of $W$ are zero can model the emergence, loss, or change of a self-loop or an edge among boundary nodes. In this case, the perturbation acts additively on $\subscr{Q}{red} = Q/Q(\alpha,\alpha)$\,as
\begin{equation}
	\subscr{\tilde Q}{red} 
	\triangleq 
	\tilde Q/\tilde Q(\alpha,\alpha)
	=
	\subscr{Q}{red} + W[\alpha,\alpha]
	\label{eq: boundary node perturbation}
	\,.
\end{equation}
If the perturbation affects only the interior nodes, then $W(\alpha,\alpha)$ is a non-zero matrix and all other elements of $W$ are zero. Inspired by \cite{BA:02,AJW-BFW:96}, we put more structure on the perturbation matrix $W$ and consider symmetric rank one perturbations of the form $W = \Delta \cdot (e_{i} - e_{j}) (e_{i} - e_{j})^{T}$, where $\Delta \in \mbb R$. Such a perturbation changes the weight of the edge $\{i,j\}$ in the network from $A_{ij}$ to $A_{ij} + \Delta$ and can model the loss or the appearance of a new link $\{i,j\}$ with weight $A_{ij} = \Delta$. Since a perturbation among the boundary nodes is fully captured by \eqref{eq: boundary node perturbation}, the more interesting case of a perturbation among the interior nodes is presented in the following theorem.

\begin{theorem}[\bf Perturbation of the Interior Network]
\label{Theorem: Effect of rank one perturbation}
Consider the Kron-reduced matrix $\subscr{Q}{red} = Q/Q(\alpha,\alpha)$, the accompanying matrix $\subscr{Q}{ac} = -Q[\alpha,\alpha)Q(\alpha,\alpha)^{-1}$, and a  symmetric rank one perturbation\footnote{The perturbation $W$ affects the edge between the $i$th and the $j$th interior node by a value $\Delta$.} $W \triangleq \Delta \cdot (e_{i+|\alpha|} - e_{j+|\alpha|}) (e_{i+|\alpha|} - e_{j+|\alpha|})^{T}$ for distinct $i,j \in \mc I_{n-|\alpha|}$ and such that the perturbed matrix $\tilde Q \triangleq Q +W$ remains an irreducible loopy Laplacian.
The following statements hold:

\begin{enumerate}

\item[1)] {\bf Algebraic perturbation:} $\subscr{Q}{red}$ undergoes the rank one perturbation
\begin{equation}
	\subscr{\tilde Q}{red} 
	\triangleq
	\tilde Q/ \tilde Q(\alpha,\alpha)
	\equiv
	\subscr{Q}{red}
	+
	\frac{\subscr{Q}{ac} (e_{i} - e_{j}) \Delta (e_{i} - e_{j})^{T} \subscr{Q}{ac}^{T}}{1 + \Delta \cdot \subscr{R}{int}[i,j]} 	
	\label{eq: interior node perturbation - final}
	\,,
\end{equation}
where $\subscr{R}{int}[i,j] \triangleq (e_{i} - e_{j})^{T} Q(\alpha,\alpha)^{-1} (e_{i} - e_{j}) \geq 0$. 

\item[2)] {\bf Resistive perturbation:} Let $R$ and $\tilde R$ be the matrices of effective resistances corresponding to $Q$ and $\tilde Q$, respectively. For any $k,\ell \in \mc I_{n}$ it holds\,\,that
\begin{equation}
	\tilde R_{k\ell} 	
	=
	R_{k\ell} - \frac{\Delta}{1+\Delta \cdot R_{i+|\alpha|,|j+ \alpha|}} \| (e_{k} - e_{\ell})^{T} Q^{\dagger} (e_{i+|\alpha|} - e_{j+|\alpha|})\|_{2}^{2}
	\label{eq: resistive perturbation bound}
	\,.
\end{equation}
If $\Delta > 0$ (respectively $\Delta<0$) then $\tilde R_{k\ell} \leq R_{k\ell}$ (respectively $\tilde R_{k\ell} \geq R_{k\ell}$).

\end{enumerate}
\end{theorem}

Note that the term $\subscr{R}{int}[i,j]$ in statement 1)  can be interpreted as the effective resistance between the perturbed nodes in the interior network 
induced by $Q(\alpha,\alpha)$. 
Likewise, the physical interpretation of the term $\subscr{Q}{ac} (e_{i}-e_{j})\Delta(e_{i}-e_{j})^{T} = \subscr{Q}{ac} W(\alpha,\alpha)$ is well-known in network theory \cite{AJW-BFW:96}. The perturbation $W$ has the same effect on\,\,the network equations $I=(Q+W)V$ as the current injection $\tilde I = -W V$, that is, the perturbation of the interior edge $\{i,j\}$ by a value $\Delta$ is equivalent to injecting the\,current $\Delta \cdot (V_{i+|\alpha|} - V_{j+|\alpha|})$  into the $j$th interior node and extracting it from the $i$th interior node. In the reduced network equations \eqref{eq: reduced network equations} the current injection $\tilde I$ translates to the current injection $\subscr{Q}{ac} \tilde I(\alpha) = -\subscr{Q}{ac} (e_{i}-e_{j})\Delta(e_{i}-e_{j})^{T} V(\alpha)$ into the boundary nodes. Elimination of $V(\alpha)$  then yields $I[\alpha] = \subscr{\tilde Q}{red} V[\alpha]$.
%
%
Finally, the additive term in identity \eqref{eq: resistive perturbation bound} resembles the  {\it sensitivity factor} in network theory \cite{ID-MP:10,AJW-BFW:96}. From Remark \ref{Remark: physical intuition about resistance}, we see that $(e_{k} - e_{\ell})^{T} Q^{\dagger} (e_{i+|\alpha|} - e_{j+|\alpha|})$ corresponds to the potential drop between nodes $k$ and $\ell$ if a unit current is injected in the $i$th interior node and extracted at the $j$th interior node. As before, the current flowing along the perturbed edge $\{i,j\}$ in the interior network is redistributed in the perturbed network in identity\,\eqref{eq: resistive perturbation bound}. 
%

Starting from identity \eqref{eq: interior node perturbation - final} various spectral bounds can be derived for the rank one perturbations of $Q$ and $\subscr{Q}{red}$. For instance, for $\Delta<0$ Weyl's inequalities \eqref{eq: Weyl's inequality}\,give
\begin{equation*}
	\lambda_{r}(\subscr{Q}{red}) 
	\geq
	\lambda_{r}(\subscr{\tilde Q}{red})
	\geq 
	\lambda_{r}(\subscr{Q}{red}) + \frac{\Delta}{1+\Delta \cdot \subscr{R}{int}[i,j]} \cdot \| \subscr{Q}{ac} (e_{i} - e_{j}) \|_{2}^{2}
	\,,\quad r \in \mc I_{|\alpha|}
	\,.
\end{equation*}
These bounds can be further related to $Q$ and $\tilde Q$ via the interlacing inequalities \eqref{eq: bound on Laplacian eigenvalues and their Schur complements} or  \cite[Theorem 4.3.4]{RAH-CRJ:85}. Finally, we remark that the proof of Theorem \ref{Theorem: Effect of rank one perturbation}, presented in the following, can also be applied to more complex symmetric perturbations. 

{\em Proof of Theorem \ref{Theorem: Effect of rank one perturbation}}.
Since the perturbed matrix $\tilde Q = Q + W$ is a symmetric and irreducible loopy Laplacian, the reduced matrix $\subscr{\tilde Q}{red} = \tilde Q/ \tilde Q(\alpha,\alpha)$ exists by Lemma \ref{Lemma: Structural Properties of Kron Reduction}. By employing the matrix identity \eqref{eq: identity for inv(A+B) - 4}, the Schur complement $\subscr{\tilde Q}{red}$ simplifies\,to
\begin{align*}
	\subscr{\tilde Q}{red} 
	&=
	\bigl( Q + W \bigr) / \bigl( Q(\alpha,\alpha) + W(\alpha,\alpha) \bigr)
	=
	\bigl( Q + W \bigr) / \bigl( Q(\alpha,\alpha) + \Delta (e_{i} - e_{j})^{T} (e_{i} - e_{j}) \bigr)
	\\
	&=
	\subscr{Q}{red}
	+
	\Delta \cdot
	\frac{Q[\alpha,\alpha) Q(\alpha,\alpha)^{-1} (e_{i}-e_{j}) (e_{i}-e_{j})^{T} Q(\alpha,\alpha)^{-1} Q(\alpha,\alpha]}{1 + \Delta \cdot (e_{i}-e_{j})^{T} Q(\alpha,\alpha)^{-1} (e_{i}-e_{j})}
	\,.
\end{align*}	
This equality further simplifies to identity \eqref{eq: interior node perturbation - final} with $\subscr{R}{int}$ as defined above. This concludes the proof of statement 1).
%
For the proof of statement 2), we initially consider the strictly loopy case. In this case, $\tilde Q^{-1}$ can be obtained from identity \eqref{eq: identity for inv(A+B) - 4}\,\,as
\begin{equation*}
	\tilde Q^{-1}
	=
	(Q+W)^{-1}
	=
	Q^{-1}
	-
	\frac{\Delta \cdot Q^{-1} (e_{i+|\alpha|} - e_{j+|\alpha|})^{T} (e_{i+|\alpha|} - e_{j+|\alpha|}) Q^{-1}}{1+\Delta (e_{i+|\alpha|} - e_{j+|\alpha|})^{T} Q^{-1} (e_{i+|\alpha|} - e_{j+|\alpha|})}
	\,.	
\end{equation*}
Multiplication of $\tilde Q^{-1}$ from the left by $(e_{k} - e_{\ell})^{T}$ and from the right by $(e_{k} - e_{\ell})$ yields  identity \eqref{eq: resistive perturbation bound}.
In the loop-less case when $Q$ is singular, the same arguments can be applied on the image of $Q$ by considering the non-singular matrix $\tilde Q + (\delta/n) \, \fvec 1_{n \times n}$\,\,for $\delta \neq 0$ and identity \eqref{eq: effective resistance from auxiliary matrix}. This results in the more general identity \eqref{eq: resistive perturbation bound}. The second part of statement 2) follows then again from Rayleigh's monotonicity law \cite{PGD-JLS:84}.
%
\qquad\endproof

\section{Conclusions}\label{Section: Conclusions}

We studied a graph-theoretic analysis of the Kron reduction process. Our
analysis is motivated by various example applications spanning from classic
circuit theory over electrical impedance tomography to power network
applications. Prompted by these applications, we presented a detailed and
comprehensive graph-theoretic analysis of the Kron reduction process. In
particular, we carried out a thorough topological, algebraic, spectral,
resistive, and sensitivity analysis of the Kron-reduced matrix. 
This analysis led to novel results in algebraic graph theory and new physical insights in the application domains of Kron reduction.
We believe our results can be directly employed in the application areas of Kron
reduction.

Of course, the results contained in this paper can and need to be further
refined to meet the specific problems and graph topologies in each
particular application area. Our analysis also demands answers to further
general questions, such as the extension of this work to directed graphs
\cite{MF:86} or complex-valued Laplacian matrices occurring in all
disciplines of circuit theory 
\cite{CSA-MHC-YII:05,BNS:07,JBW:09,PWS-MAP:98}. The results on the
effective resistance should be extended to the related but more general
network sensitivity factors \cite{ID-MP:10,AJW-BFW:96}. The presented topological
analysis should be rendered to a synthesis of optimal reduction patterns in
terms of sparsity similar to
\cite{AR:96,JR-WHAS:09,MF:76,JRG:94,YS:03,YS-MS:00}. Some of our results have
also analogous results for the Perron complement \cite{CDM:89-2,MN:00}
occurring in Markov chain applications \cite{CDM:89}, but many results
still have to be mapped from Perron complements to Schur complements and
vice versa. Finally, it would be of interest to analyze the effects of Kron
reduction on centrality measures, clustering coefficients, and more general
graph-theoretic metrics than the effective resistance.




\bibliographystyle{siam}
\bibliography{alias,Main,FB,New}


\end{document}